\documentclass[12pt]{article}
 \usepackage[english]{babel}
 \usepackage{amsmath}
 \usepackage{amssymb}
 \usepackage{amsthm}
  \usepackage{amsmath}
\usepackage{amsbsy}
\usepackage{amsthm}
\usepackage{amssymb}
\usepackage{enumerate}
\usepackage{xspace}
\usepackage{euscript}
\usepackage{graphicx}
\usepackage{amscd}
 
 \newtheorem{lemma}{Lemma}[section]
 \newtheorem{proposition}[lemma]{Proposition}
 \newtheorem{definition}[lemma]{Definition}
 \newtheorem{theorem}[lemma]{Theorem}
 \newtheorem{corollary}[lemma]{Corollary}
 \newtheorem{remark}[lemma]{Remark}

 \newcommand{\no}{{\nonumber}}
 \newcommand{\bear}{\begin{array}}
 \newcommand{\enar}{\end{array}}
 \newcommand{\beq}{\begin{equation}}
 \newcommand{\eeq}{\end{equation}}
 \newcommand{\beqn}{\begin{eqnarray}}
 \newcommand{\eeqn}{\end{eqnarray}}
 \newcommand{\beit}{\begin{itemize}}
 \newcommand{\eeit}{\end{itemize}}
 \newcommand{\beal}{\begin{align}}
 \newcommand{\eeal}{\end{align}}

 \newcommand{\rsp}{{\bf R}}
 \newcommand{\csp}{{\bf C}}
 \newcommand{\nsp}{{\bf N}}

 \newcommand{\ds}{\displaystyle}
 \newcommand{\ve}{\varepsilon}
 \newcommand{\g}{\gamma}
 \renewcommand{\l}{\lambda}
 \newcommand{\dl}{\delta}
 
 \newcommand{\s}{\sigma}
 \newcommand{\G}{\Gamma}
 \newcommand{\om}{\omega}
 \newcommand{\Om}{\Omega}
 \renewcommand{\a}{\alpha}
 \renewcommand{\b}{\beta}

 \newcommand{\ov}{\overline}

 \newcommand{\wtil}[1]{\widetilde{#1}}
 \newcommand{\pn}{\par \noindent}
 \newcommand{\med}{\medskip}
 \newcommand{\qq}{\qquad}
 \newcommand{\q}{\quad}

 \newcommand{\de}{\,{\rm d}}

\title{{Generation type inequalities\\
for closed linear operators\\
related to domains with conical points}}

\author{Alberto Favaron (Bologna)
\footnote{The author is member
of GNAMPA of the Italian Istituto Nazionale di Alta Matematica
(INdAM).}}

\date{}

\begin{document}
\maketitle \pn

{\bf Abstract.}
Let ${\cal A}(x;D_x)$ be a second-order linear differential operator
in divergence form.
We prove that the operator $\l I- {\cal A}(x;D_x)$, 
where $\l\in\csp$ and $I$ stands for the identity operator, 
is closed and injective when ${\rm Re}\l$ is large enough and the domain
of ${\cal A}(x;D_x)$
consists of a special class of weighted Sobolev
function spaces related to conical open bounded sets of $\rsp^n$, 
$n \ge 1$.
\med \pn
{\it Key words and phrases.} Resolvent estimates.
Weighted Sobolev function spaces.
Conical bounded domains of $\rsp^n$.
\section{Introduction and plan of the paper}
\setcounter{equation}{0} 
In this paper we present a new approach
for proving an estimate of generation type for the norm of the
resolvent $[\l I- {\cal A}(x;D_x)]^{-1}$ of the operator
$\l I- {\cal A}(x;D_x)$,  where $\l\in\csp$, $I$ stands for the identity
operator and ${\cal A}(x;D_x)$ denotes 
the second-order linear differential operator in divergence form
\beqn
\label{1.1} {\cal A}(x;D_x) =\sum_{j=1}^{n}D_{x_j}
\big(\sum_{k=1}^{n}a_{j,k}(x)D_{x_k} \big).
\eeqn
We stress that in our paper the domain of ${\cal A}(x;D_x)$
will consist of an appropriate class of weighted
Sobolev spaces whose elements will be functions taking 
their values in conical open bounded sets of $\rsp^n$, $n \ge 1$.
\\
With the language of the modern semigroup 
theory  a generation type estimate means that,
denoted with ${\cal L}(X)$ the Banach space of the linear
bounded operators from $X$ to $X$, $X$ being a Banach space,
and endowed ${\cal L}(X)$ with the usual uniform
operatorial norm, then
 $\|[\l I- {\cal A}(x;D_x)]^{-1}\|_{{\cal L}(X)}$
is bounded from above
by some constant times $|\l|^{-1}$, at least
for large enough ${\rm Re}\l$.
\\
Even if in order to prove our main result we adopt
an idea that goes back to \cite{A} and \cite{AN},
i.e. the procedure of increasing the dimension from $n$ to $n+1$,
in our proof there are so many different elements
with respect to the proof of the estimates
in the quoted papers that we may consider
our results totally independent of those.
\\
The novelties arise fundamentally from the fact
that we consider bounded domain
$G$ of $\rsp^n$ having a singular point $O$ situated
in a part of the boundary $\partial G$ with a conical
structure. This forces us to consider weighted Sobolev
function spaces for which, unfortunately, the classical {\it a priori}
estimates of \cite{ADN1} are not available.
The role of that estimates will be played here
by some estimates of the same type proven in \cite{MP3},
but these estimates, when applied
in dimension $n+1$, require
the conical structure of $G$ to be preserved
when we increase the dimension.
Therefore, denoted by $x_0$ the added variable,
unlike standard procedures making use of
 $\Gamma=(-\infty,+\infty)\times G$, 
we will consider, as a new domain
in dimension $n+1$, a domain $\wtil G$ which can be regarded 
as a rotation of $G$ around its symmetry axis.
\\
Of course, when $n=1$, $G$ simply coincides with a bounded
open interval of $\rsp$ and rotations have no meaning.
However, in this situation too,
it will always be possible to consider
two-dimensional conical domains $\wtil G$ having $G$ as their symmetry axis.
\\
The main difficulties arising by the use of $\wtil G$
instead of that of $\Gamma$ consist in the following:
\begin{itemize}
\item[(i)] 
the proof that the property of the boundary conditions
to cover ${\cal A}(x;D_x)$ on $\partial G\backslash\{O\}$
in the sense of \cite{ADN2} continues to hold when we increase the
dimension. This is not a straightforward fact and
forces us to implement a new set of boundary operators coinciding
with the original ones on $\partial G$;
\item[(ii)] 
the necessity of considering cut-off functions depending on both
 variables $x_0$, $x$, where $x\in G$,
instead of cut-off functions as those considered 
in \cite{A} and \cite{AN} and depending 
only on the added variable $x_0$.
As a consequence, our computations will be heavier and longer than those
in the quoted papers (cf. also \cite{LU}).
\end{itemize}
Observe that we will consider bounded domains of conical type.
Since a lot of papers have been devoted
in the past to the investigation of boundary value problems
in such domains, we
prefer here not to mention any of them, but only to refer the
interested reader to \cite{MP1},
where some examples of admissible domains
and an exhaustive list of references for
this kind of problems are given.
\\
We would like to emphasize that
generation type estimates are one of the
main tools needed to prove that a linear operator generates an
analytic semigroup of linear bounded operators. Hence,
if the generation in our functional setting
could be guaranteed, by showing the
surjectivity of $\l I- {\cal A}(x;D_x)$ too,
the range of applications of our result would be extremely large.
Indeed, nowadays semigroup theory is one of the most used tool
in both direct and inverse problems
related parabolic differential equations.
However, while for regular
domains and classical Sobolev spaces many generation
results are available, the same, to the author's knowledge,
is not true for conical domains and weighted Sobolev spaces.
\\
The plan of this paper is the following.
In Section 2, using notations of \cite{MP3},
we introduce the class of domains
and of weighted Sobolev function spaces
we will deal with.
Moreover, we introduce also the
correspondent spaces of traces for the boundary values.
\\
Section 3 is devoted to recall the {\it a priori} estimates of \cite{MP3}
for boundary value problems in the functional setting
of Section 2. To this purpose we need to introduce
some further technical definitions and a rather heavy notation 
which, however, having to deal with scalar and not matrix
differential operators,
turns out to be quite simple in our case.
\\
In the first part of Section 4
we list all the basic assumptions on the domain $G$,
on the operator ${\cal A}(x;D_x)$ and on the
boundary operator ${\cal B}(x;D_x)$ associated with ${\cal A}(x;D_x)$.
Under these assumptions, in the second part of Section 4
we will introduce the concept of regular boundary value problem
and, for such a problem, we will state our main result (Theorem 4.6).
We conclude the section by showing some easy corollaries
to our estimate and related to the analytic semigroups theory. 
\\
Section 5 contains the proof of the preliminary Lemma \ref{lem4.2}.
Essentially, it states that the property of ${\cal B}(x;D_x)$
to cover ${\cal A}(x;D_x)$ on $\partial G\backslash\{O\}$
in the sense of \cite{ADN2} continues to hold when we increase the
dimension, provided we replace 
the triplet $\{{\cal A}(x;D_x),{\cal B}(x;D_x),G\}$
with the triplet $\{{\cal A}(x;D_x)+{\rm e}^{i\psi}D_{x_0}^2,
{\cal B}(x;D_x)+x_0D_{x_0},\wtil G\}$, $\psi\in[-\pi/2,\pi/2]$.
\\
In Section 6 we introduce the class of our 
admissible cut-off functions.
For the reasons we said before, they have a structure 
more complicated (cf. (\ref{6.2})) than those used in \cite{A} and \cite{AN} 
and hence, for clarity's sake, we report 
all the necessary computations we need in order to perform 
the technicalities of Section 7.
\\
Finally, in Section 7 we prove our main result.
The proof will be derived simply by
taking advantage of the assumptions on $G$ and
by combining Theorem \ref{3.2}
with Lemma \ref{lem4.2} and with the further preliminary
estimates of Lemma \ref{lem7.3} and Lemma \ref{lem7.4}.
\section{The spaces $V_{p,\b}^l(G)$, $W_{p,\b}^l(G)$,
$V_{p,\b}^{l-p^{-1}}(\partial G)$, $W_{p,\b}^{l-p^{-1}}(\partial G)$
}
\setcounter{equation}{0}
Let $B(0,1)$ be the unit open ball of $\rsp^n$, $n\ge 1$,
and denote by $K$ an open cone of $\rsp^n$ having its vertex
at the origin and cutting out on the unit sphere
$\partial\mbox{}B(0,1)$ a domain $\Om$.\\
From now on, with $G$ we will denote an open subset of $\rsp^n$
having compact closure $\ov G$ and boundary $\partial\mbox{}G$
on which there is a point $O$ such that:
\\[1mm]
\indent
(i)\ \ $\partial\mbox{}G\backslash\{O\}$ is a smooth, 
$(n-1)$--dimensional submanifold of $\rsp^n$;
\\[1mm]
\indent
(ii)\ \ near $O$ the domain $G$ coincides with $K\cap B(0,1)$.
\\[2mm]
Using a multi-index notation, for $1<p<+\infty$, $\b\in\rsp$, $l=0,1,\ldots$,
we define the weighted spaces $V_{p,\b}^l(G)$ and $W_{p,\b}^l(G)$
as the spaces of functions $u$ in $G$ endowed, respectively,
with the following
norm $\|\cdot\|_{V_{p,\b}^l(G)}$ and $\|\cdot\|_{W_{p,\b}^l(G)}$,
where $|x|=(x_1^2+\ldots+x_n^2)^{1/2}$:
\beqn\label{2.1}
\hskip -1truecm
\|u\|_{ V_{p,\b}^l(G)}
\!\!\!&=&\!\!\!
\Big(\sum_{\a=0}^l\int_{G}
|x|^{p(\b-l+|\a|)}|D^{\a}u(x)|^p\de x\Big)^{1/p}<+\infty
\\[2mm]
\label{2.2}
\hskip -1truecm
\|u\|_{ W_{p,\b}^l(G)}
\!\!\!&=&\!\!\!
\Big(\sum_{\a=0}^l\int_{G}
|x|^{p\b}|D^{\a}u(x)|^p\de x\Big)^{1/p}<+\infty
\eeqn
Since (\ref{2.1}) and (\ref{2.2}) coincide if  $l=0$ we set
$L_{p,\b}(G)$ to be the weighted $L_p$ space of functions in $G$
endowed with norm
\beqn\no
\hskip -1truecm
\|u\|_{ L_{p,\b}(G)}
\!\!\!&=&\!\!\!
\Big(\int_{G}
|x|^{p\b}|u(x)|^p\de x\Big)^{1/p}\,.
\eeqn
As shown in \cite{MP3}, the space
$C_0^{\infty}(\ov G\backslash\{O\})$ of the infinitely differentiable
functions having compact support on $\ov G\backslash\{O\}$
is dense in $V_{p,\b}^l(G)$ and the following theorem holds true.
\begin{theorem}\label{thm2.1}
1) If $\b<-np^{-1}$ or $\b>l-np^{-1}$
then the spaces $V_{p,\b}^l(G)$ and $W_{p,\b}^l(G)$ coincide
and the norm (\ref{2.1}), (\ref{2.2}) are equivalent.\\
2) If for some number $\nu=0,1,\ldots,l-1$ the inequalities
$\nu-np^{-1}<\b<\nu+1-np^{-1}$ are satisfied, then the space
$W_{p,\b}^l(G)$ is the direct sum of $V_{p,\b}^l(G)$
and $\Pi_{l-\nu-1}$, where $\Pi_{l-\nu-1}$ is
the space of polynomials in $x$ of degree
at most equal to $l-\nu-1$ $(\Pi_q=\{0\}$ if $q<0)$.
\end{theorem}
\begin{proof}
See the proof of Theorem 2.1 in \cite{MP3}.
\end{proof}
In order to consider boundary value problems
we need to define also the spaces $V_{p,\b}^{l-p^{-1}}(\partial\mbox{}G)$
and  $W_{p,\b}^{l-p^{-1}}(\partial\mbox{}G)$, i.e. the
spaces of traces on $\partial\mbox{}G$ of functions in $V_{p,\b}^l(G)$
and $W_{p,\b}^l(G)$, respectively. It turns out that
$V_{p,\b}^{l-p^{-1}}(\partial\mbox{}G)$ is the
quotient space $V_{p,\b}^l(G)\backslash\wtil V_{p,\b}^l(G)$, where
$\wtil V_{p,\b}^l(G)$ is the completion
with respect to the $V_{p,\b}^l(G)$--norm
of the set of smooth functions in $V_{p,\b}^l(G)$ equal to
zero on $\partial\mbox{}G$.
$V_{p,\b}^{l-p^{-1}}(\partial\mbox{}G)$ is endowed with the
norm:
\beqn\label{2.3}
\|u\|_{V_{p,\b}^{l-p^{-1}}(\partial\mbox{}G)}\,=\,
\inf\big\{\|v\|_{V_{p,\b}^l(G)}:v-u\in \wtil V_{p,\b}^l(G)\big\}\,.
\eeqn
Replacing $V$ with $W$ in the above definitions
we obtain the description of $W_{p,\b}^{l-p^{-1}}(\partial\mbox{}G)$.\\
From the fact that $C_0^{\infty}(\ov G\backslash\{O\})$
is dense in $V_{p,\b}^l(G)$ it easily follows that
$C_0^{\infty}(\partial\mbox{}\ov G\backslash\{O\})$
is dense in $V_{p,\b}^{l-p^{-1}}(\partial\mbox{}G)$.
Moreover, Theorem \ref{thm2.1} ensures that
if $\b<-np^{-1}$ or $\b>l-np^{-1}$ then
$V_{p,\b}^{l-p^{-1}}(\partial\mbox{}G)$ and
$W_{p,\b}^{l-p^{-1}}(\partial\mbox{}G)$ coincide whereas
(cf. \cite[Theorem 3.1]{MP3})
if for some number $\nu$ the inequalities $\nu-np^{-1}<\b<\nu+1-np^{-1}$
are satisfied then $W_{p,\b}^{l-p^{-1}}(\partial\mbox{}G)$
is the direct sum of $V_{p,\b}^{l-p^{-1}}(\partial\mbox{}G)$
and the space $Y_{l-\nu-1}$ of polynomials of degree at most $l-\nu-1$
which are not identically zero on
$\partial\mbox{}\Om\times\rsp_+$.
\section{Admissible operators and\\ boundary value problems in
$W_{p,\b}^l(G)$}
\setcounter{equation}{0}
Let  $D_x$ denotes the $n$-uple $(D_{x_1},\ldots,D_{x_n})$
and let ${\mathfrak C}(\mu,s)$ to be the class
of differential operators ${\cal M}(x; D_x)$
of order $\mu$ with coefficients in $C^s(\ov G\backslash\{O\};\csp)$
and admitting, near $O$, the following representation
in local spherical co-ordinates $(r,\om)$:
\beqn\label{3.1}
&&\hskip -1truecm
{\cal M}(x, D_x)\,
=\,r^{-\mu}\sum_{k+|\g|\le\mu}
p_{k,\g}(r,\om)(rD_r)^kD_\om^\g
\,\equiv\,
r^{-\mu}M(r,\om;rD_r,D_\om)\,,
\eeqn
where the functions $p_{h,\a}(r,\om)$, $h+|\a|\le\mu$, satisfy
the condition
\beqn\label{3.2}
(rD_r)^qD_\om^\g\,p_{h,\a}\in C([0,\dl]\times\ov\Om;\csp)\,,\q\;
q+|\g|\le s,\ \dl={\rm const}>0.
\eeqn
Recall that, for $h>0$, we have
\beqn\label{3.3}
&&\hskip -1,5truecm
D_r^h\,=\sum_{|\a|=h}a_{h,\a}(\om)D_x^{\a}\,,\qq
D_{\om}^h\,=\sum_{0<|\a|\le h}r^{|\a|}b_{h,\a}(\om)D_x^{\a}\,,
\eeqn
where $a_{h,\a}$ and $b_{h,\a}$, $0<|\a|\le h$,
are smooth functions on $\partial\mbox{}B(0,1)$.\\
From (\ref{2.1}) it is easy to prove that
any ${\cal M}\in{\mathfrak C}(\mu,s)$
realizes a continuous mapping $V_{p,\b}^l(G)\to V_{p,\b}^{l-\mu}(G)$
for $s\ge l-\mu$ and, by Theorem \ref{thm2.1},
if $\b<-np^{-1}$ or $\b>l-np^{-1}$ the same property holds true
with $V_{p,\b}^l(G)$ and $V_{p,\b}^{l-\mu}(G)$ replaced
by $W_{p,\b}^l(G)$ and $W_{p,\b}^{l-\mu}(G)$, respectively. In the
case there exists $\nu=0,\ldots,l-1$ such that $\nu-np^{-1}<\b<\nu+1-np^{-1}$
the map ${\cal M}:W_{p,\b}^l(G)\to W_{p,\b}^{l-\mu}(G)$ is still continuous
if ${\cal M}(\Pi_{l-\nu-1})\subset W_{p,\b}^{l-\mu}(G)$.
\begin{remark}\label{rem3.1}
\emph{For reasons that will be clearer in Section 4,
denote by $(x_0,x)$
the points of $\rsp^{n+1}$ and by
$\wtil G$ an $(n+1)$-dimensional domain which,
close to the origin, coincides with the cone
$\{(x_0,x)\in\rsp^{n+1}:|(x_0,x')|\le C_1x_n,\ 
x'=(x_1,\ldots,x_{n-1}),\ x_n>0\}$, $C_1>0$.
We will show here that if $a_{j,k}\in C^1(\ov G;\csp)$, $j,k=1,\ldots, n$,
and $G=\{(x_0,x)\in \wtil G:x_0=0\}$,
then the
operator ${\cal A}_{\psi}(x;D_x, D_{x_0})=
{\cal A}(x;D_x)+e^{i\psi}D_{x_0}^2$,
where  $\psi\in[-\pi/2,\pi/2]$ and 
${\cal A}(x;D_x)$ is defined by (\ref{1.1}), 
belongs to the class ${\frak C}(2,0)$.
\\
We introduce in the space $\rsp^{n+1}$ the $(n+1)$-dimensional
spherical co-ordinates, related to the Cartesian ones by
the well-known relationships:
\beqn\label{3.4}
&&\hskip -1,5truecm
\begin{array}{lll}
&&\!\!
(x_0,x_1,\ldots,x_{n-1},x_n)
\\[2mm]
&=&\!\!
\big(r\cos\theta_0,r\sin\theta_0\cos\theta_1,\ldots,
r\prod_{h=0}^{n-2}\sin\theta_h\cos\theta_{n-1},
r\prod_{h=0}^{n-1}\sin\theta_h\big)\,,
\end{array}
\eeqn
where $r=|(x_0,x)|$, $\theta_h\in [0,\pi]$, $h=0,\ldots,n-2$,
$\theta_{n-1}\in [0,2\pi)$ and where
$\prod_{h=l_1}^{l_2}\sin\theta_h$ has to be intending equal to one
if $l_2<l_1$.\\
Denoting by $\om$ the $(n-1)$-uple
$(\theta_0,\ldots,\theta_{n-1})$,
with the help of (\ref{3.4}) it is not too difficult to show that
the gradient $(D_{x_0},D_x)$ can be expressed in
terms of $(D_r,D_{\om})$, $D_{\om}=(D_{\theta_0},\ldots,D_{\theta_{n-1}})$,
by the following formulae:
\beqn \label{3.5}
\begin{array}{lll}
&&\hskip -1truecm
D_{x_j}=\tau_{j,r}(\om)D_r
+r^{-1}\big(\sum_{k=0}^{n-2}\tau_{j,\theta_k}(\om)D_{\theta_k}
+\tau_{j,\theta_{n-1}}(\om)D_{\theta_{n-1}}\big)\,,\q\; j=0,\ldots,n\,,
\end{array}
\eeqn
where, $\dl_{i,l}$ standing for the Kronecher symbol,
for any $j\in\{0,\ldots,n-1\}$ and any
$k\in\{0,\ldots,n-2\}$ we have
\beqn\label{3.6}
&&\hskip -0,7truecm
\left\{\!\!
\begin{array}{lll}
\tau_{j,r}(\om)=\prod_{h=0}^{j-1}\sin\theta_h\cos\theta_j\,,
\q\,
\tau_{n,r}(\om)=\prod_{h=0}^{n-1}\sin\theta_h\,,
\\[3mm]
\tau_{j,\theta_k}(\om)=\left\{
\begin{array}{lll}
&\hskip -0,4truecm
\big(\prod_{h=0}^{k}\sin\theta_h\big)^{-1}
\big[\prod_{h=k}^{j-1}\sin\theta_h\cos\theta_k\cos\theta_j-\dl_{k,j}\big]
\,,\q\;{\rm if}\;\;k\le j\,,
\\[3mm]
&\hskip -0,4truecm
0\,,\q\;{\rm if}\;\;k>j\,,
\end{array}
\right.
\\[6mm]
\tau_{n,\theta_k}(\om)=
\big(\prod_{h=0}^{k-1}\sin\theta_h\big)^{-1}
\prod_{h=k+1}^{n-1}\sin\theta_h\cos\theta_k\,,
\\[3mm]
\tau_{j,\theta_{n-1}}(\om)=
-\big(\prod_{h=0}^{n-2}\sin\theta_h\big)^{-1}\dl_{(n-1),j}\sin\theta_{j}\,,
\q \,
\tau_{n,\theta_{n-1}}(\om)=
\big(\prod_{h=0}^{n-2}\sin\theta_h\big)^{-1}\cos\theta_{n-1}\,.
\end{array}
\right.
\eeqn
Hence, if we set $\wtil a_{j,k}(r,\om)=a_{j,k}(r\sin\theta_0\cos\theta_1,
\ldots,r\prod_{h=0}^{n-2}\sin\theta_{h}\cos\theta_{n-1},
r\prod_{h=0}^{n-1}\sin\theta_h)$,
$j,k=1,\ldots,n$, we obtain
\beqn\label{3.7}
\begin{array}{lll}
&&\hskip -1,5truecm
\sum_{k=1}^{n}a_{j,k}(x)D_{x_k}=
f_{j,r}(r,\om)D_r
+r^{-1}\sum_{h=0}^{n-1}
f_{j,\theta_h}(r,\om)D_{\theta_h},\q\; j=1,\ldots,n\,,
\end{array}
\eeqn
functions $f_{j,r}$, $f_{j,\theta_l}$, $j=1,\ldots,n$, $l=0,\ldots,n-1$,
being defined by
\beqn\label{3.8}
\begin{array}{lll}
&&\hskip -1,5truecm
f_{j,r}(r,\om):=
\sum_{k=1}^n{\wtil a}_{j,k}(r,\om)\tau_{k,r}(\om)\,,
\qq
f_{j,\theta_l}(r,\om):=
\sum_{k=1}^n{\wtil a}_{j,k}(r,\om)\tau_{k,\theta_l}(\om)\,.
\end{array}
\eeqn
Using again (\ref{3.5}) and applying it to relations
(\ref{3.7}), performing easy computations we get
\beqn\label{3.9}
\begin{array}{lll}
&&\hskip -1,5truecm
{\cal A}(x;D_x)={\cal Q}_r(r,\om;D_r,D_{\om})
+\sum_{l=0}^{n-1}{\cal Q}_{\theta_l}(r,\om;D_r,D_{\om}),
\end{array}
\eeqn
where ${\cal Q}_r(r,\om;D_r,D_{\om})$,
${\cal Q}_{\theta_l}(r,\om;D_r,D_{\om})$, $l=0,\ldots,n-1$, stand
for the second-order linear differential operator
\beqn
\label{3.10}
\begin{array}{lll}
&&\hskip -1,3truecm
{\cal Q}_r(r,\om;D_r,D_{\om})=
D_r\big(k_r(r,\om)D_r
+r^{-1}\sum_{h=0}^{n-1}k_{\theta_h}(r,\om)D_{\theta_h}\big),
\\[2mm]
&&\hskip -1,3truecm
{\cal Q}_{\theta_l}(r,\om;D_r,D_{\om})=r^{-1}\sum_{j=1}^n
\tau_{j,\theta_l}(\om)D_{\theta_l}
\big(f_{j,r}(r,\om)D_r
+r^{-1}\sum_{h=0}^{n-1}f_{j,\theta_{h}}(r,\om)D_{\theta_h}\big),
\end{array}
\eeqn
functions $k_r$, $k_{\theta_l}$, $l=0,\ldots,n-1$,
appearing in (\ref{3.10}) being defined by
\beqn\label{3.12}
\begin{array}{lll}
&&\hskip -1,5truecm
k_r(r,\om):=
\sum_{j=1}^n\tau_{j,r}(\om)f_{j,r}(r,\om)\,,
\qq
k_{\theta_l}(r,\om):=
\sum_{j=1}^n\tau_{j,r}(\om)f_{j,\theta_l}(r,\om)\,.
\end{array}
\eeqn
Moreover, since (\ref{3.5}), (\ref{3.6}) imply
$D_{x_0}=\cos\theta_0D_r-r^{-1}{\sin\theta_0}D_{\theta_0}$,
taking advantage from
\beqn\label{3.13}
(rD_r)^2=r^2D_r^2+rD_r
\eeqn
we obtain
\beqn\label{3.14}
\begin{array}{lll}
&&\hskip -1truecm
D_{x_0}^2=r^{-2}\{\cos^2\theta_0(rD_r)^2
+\sin(2\theta_0)[I-(rD_r)]D_{\theta_0}+\sin^2\theta_0 D_{\theta_0}^2
-\cos(2\theta_0)(rD_r)\}.
\end{array}
\eeqn
Therefore, assuming $a_{i,j}\in C^1(\ov G;\csp)$, $i,j=1,\ldots,n$,
taking into account the following formulae (cf. (\ref{3.3}) with $h=1$)
\beqn\label{3.15}
&&\hskip -1,5truecm
\left\{\!\!
\begin{array}{lll}
D_r\!\!&=&\!\!
\sum_{k=0}^{n-1}(\prod_{h=0}^{k-1}\sin\theta_h\cos\theta_k) D_{x_k}
+\prod_{h=0}^{n-1}\sin\theta_h D_{x_n}\,,
\\[2mm]
D_{\theta_l}\!\!&=&\!\!\sum_{j=l}^n(D_{\theta_k}x_j)D_{x_j}\,,
\qq l=0,\ldots,n-1\,,
\end{array}\right.
\eeqn
and recalling (\ref{3.8}) and (\ref{3.12}),
if we differentiate, respectively with respect to $r$
and $\theta_l$, $l=0,\ldots,n-1$, each term in the brackets of
(\ref{3.10}) and we rearrange the term using (\ref{3.13}),
from (\ref{3.9}), (\ref{3.14}) we can easily see that
${\cal A}_{\psi}(x;D_x,D_{x_0})={\cal A}(x;D_x)+e^{i\psi}D_{x_0}^2$,
admits representation (\ref{3.1}) with $\mu=2$.
In addition, since the points
$(0,\ldots,0,x_j,0,\ldots,0)$, $j=0,\ldots,n-1$, with $x_j\neq 0$
do not belong to $\wtil G$, we have
$\sin\theta_j\neq 0$ for any $j\in\{0,\ldots,n-1\}$ and hence
(cf. (\ref{3.6})) condition (\ref{3.2}) is satisfied, too, with $s=0$.}
\end{remark}
Coming back to our purposes, we consider the
boundary value problem:
\beqn\label{3.16}
&&\hskip -1truecm
{\cal L}(x;D_x)u\,=\,{\cal F}\,\q {\rm in}\;\, G\;;
\qq
{\cal B}(x;D_x)u\,=\,{\cal G}\,
\q {\rm on}\;\, \partial\mbox{}G\backslash\{O\}\,,
\eeqn
${\cal L}$ and ${\cal B}$ being matrix differential operators in $G$
of dimension $k\times k$ and $m\times k$, respectively, with elements
${\cal L}_{h,j}(x;D_x)$ and
${\cal B}_{q,j}(x;D_x)$,
$h,j=1,\ldots,k$, $q=1,\ldots,m$.
The orders of operators ${\cal L}_{h,j}$ and ${\cal B}_{q,j}$
are equal to $(s_h+t_j)$ and $(\s_q+t_j)$, respectively,
where $\{s_h\}_{h=1}^k$, $\{t_j\}_{j=1}^k$ and $\{\s_q\}_{q=1}^m$
are collections of integers with $\max_{h=1,\ldots,k}s_h=0$,
$t_j>0$, $j=1,\ldots,k$, and
$\sum_{j=1}^k(s_j+t_j)=2m$.
Clearly, ${\cal L}_{h,j}\equiv 0$ and ${\cal B}_{q,j}\equiv 0$
if $s_h+t_j<0$ and $\s_q+t_j<0$.
Moreover, taken $l\ge\max\{0,\max_{q=1,\ldots,m}\s_q\}$,
we assume ${\cal L}_{h,j}(x;D_x)\in{\mathfrak C}(s_h+t_j, l-s_h)$
and ${\cal B}_{q,j}(x;D_x)\in{\mathfrak C}(\s_q+t_j, l-\s_q)$
in a neighborhood of $O$ (cf. \cite{MP2} p.76).
\\
We require ${\cal L}$ to be {\it uniformly elliptic in
$\ov G\backslash\{O\}$} in
the sense of \cite{ADN2} and we impose that
the boundary conditions ${\cal B}$
cover ${\cal L}$ on $\partial\mbox{}G\backslash\{O\}$
(cf. \cite{ADN2} or \cite{LM}).\\
Problem (\ref{3.16}) generates a model problem in the cone
$(0,+\infty)\times \Om$. With the pair $\{{\cal L}, {\cal B}\}$
we associate the operator
${\cal U}(0,\om,z, D_{\om})=\{L(0,\om,z, D_{\om}),B(0,\om,z, D_{\om})\}$
where $\om\in\Om$, $z\in\csp$
and the matrix differential operators $L(0,\om,z, D_{\om})$
and $B(0,\om,z, D_{\om})$ are defined, respectively, by
\beqn
\label{3.17}
\hskip -1truecm
L(0,\om,z, D_{\om})\!\!&=&\!\!
\big(L_{h,j}(0,\om;z-it_j,D_{\om})\big)_{h=1.\ldots,k}^{j=1.\ldots,k}\,,
\\[2mm]
\label{3.18}
B(0,\om,z, D_{\om})\!\!&=&\!\!
\big(B_{q,j}(0,\om;z-it_j,D_{\om})\big)_{q=1,\ldots,m}^{j=1,\ldots,k}\,,
\eeqn
the operators $L_{h,j}$ and $B_{q,j}$ being determined
from ${\cal L}_{h,j}$ and ${\cal B}_{q,j}$ by means of
(\ref{3.1}) replacing $p_{k,\g}(r,\om)$ with
$p_{k,\g}(0,\om)$. As shown in the Appendix the ellipticity
of system (\ref{3.16}) implies that ${\cal U}(0,\om,z, D_{\om})$ is elliptic
with parameter in the sense of \cite{AV}.\newpage\pn
Denoting by $\vec t$, $\vec s$ and $\vec\s$ the vectors
$(t_1,\ldots,t_k)$, $(s_1,\ldots,s_k)$ and $(\s_1,\ldots,\s_m)$, 
respectively, we introduce the spaces of vector-valued functions
\beqn\label{3.19}
&&\hskip -2truecm
\begin{array}{ccc}
&V_{p,\b}^{l+\vec t}(G)=\prod_{j=1}^k V_{p,\b}^{l+t_j}(G)\,,\qq
V_{p,\b}^{l-\vec s}(G)=\prod_{j=1}^k V_{p,\b}^{l-s_j}(G)\,,&
\\[2mm]
&V_{p,\b}^{l-\vec\s-p^{-1}}(\partial\mbox{}G)=
\prod_{q=1}^m V_{p,\b}^{l-\s_q-p^{-1}}(\partial\mbox{}G)\,,&
\end{array}
\eeqn
and the correspondent spaces
$W_{p,\b}^{l+\vec t}(G)$,
$W_{p,\b}^{l-\vec\s-p^{-1}}(\partial\mbox{}G)$
obtained by replacing $V$ with $W$ in (\ref{3.19}).
By the assumptions
it is obvious that the map
\beqn\label{3.20}
&&\hskip -1truecm
\{{\cal L},{\cal B}\}:
V_{p,\b}^{l+\vec t}(G)\to V_{p,\b}^{l-\vec s}(G)
\times V_{p,\b}^{l-\vec\s-p^{-1}}(\partial\mbox{}G)\,,
\eeqn
is continuous and from Theorem \ref{thm2.1} we deduce that if
$\b<-np^{-1}$ or $\b>l+t_{\max}-np^{-1}$,
$t_{\max}=\max_{j=1,\ldots,k}t_j$, the same
regularity holds true by replacing $V$ with $W$.
In addition, Theorem 4.2 in \cite{MP3} shows that if there exists
some $\nu=0,1,\ldots,l+t_{\max}-1$ such that
$\nu-np^{-1}<\b<\nu+1-np^{-1}$ then the map
$\{{\cal L},{\cal B}\}:W_{p,\b}^{l+\vec t}(G)\to W_{p,\b}^{l-\vec s}(G)
\times W_{p,\b}^{l-\vec\s-p^{-1}}(\partial\mbox{}G)$ still remains continuous
provided that $\{{\cal L},{\cal B}\}(\Pi_{l+\vec t-\nu-1})\subset
W_{p,\b}^{l-\vec s}(G)\times
W_{p,\b}^{l-\vec\s-p^{-1}}(\partial\mbox{}G)$,
where $\Pi_{l+\vec t-\nu-1}=\prod_{j=1}^k\Pi_{l+t_j-\nu-1}$.
We can now state the following result corresponding to
Theorem 4.3 in \cite{MP3} and to which we refer the reader for the proof.
\begin{theorem}\label{thm3.2}
If $\b\notin[-np^{-1}, l+t_{\max}-np^{-1}]$
or if there exists $\nu=0,1,\ldots,l+t_{\max}-1$
such that $\nu-np^{-1}<\b<\nu+1-np^{-1}$
then for $p\in (1,+\infty)$ the operator (\ref{3.20})
is Fredholm if and only if the line ${\rm Im}z=\b-l+np^{-1}$
contains no poles of the operator ${\cal U}(0,\om,z, D_{\om})^{-1}$
which are the eigenvalues of ${\cal U}(0,\om,z, D_{\om})$.
Under this condition, for any vector-valued function
$w\in W_{p,\b}^{l+\vec t}(G)$
the following estimate holds
\beqn
\label{3.21}
\|w\|_{W_{p,\b}^{l+\vec t}(G)}
\,\le\,
c_1\big\{
\|{\cal L}w\|_{W_{p,\b}^{l-\vec s}(G)}
+\|{\cal B}w\|_{W_{p,\b}^{l-\vec\s-p^{-1}}(\partial\mbox{}G)}
+\|w\|_{W_{p,\b}^{l-1+\vec t}(G)}\big\}.
\eeqn
\end{theorem}
In the next, Theorem \ref{thm3.2} will be applied to the case in which
${\cal L}$ and ${\cal B}$ are single and not matrix differential operators.
Therefore, from now on the parameters appearing
from formula (\ref{3.16}) onward will be assumed to be the following:
\beqn\label{3.22}
k=m=1\,,\q s_1=0\,,\q t_1=2
\q \s_1=-1\,,\q l=0\,.
\eeqn
\section{Basic assumptions and main result}
\setcounter{equation}{0}
With all the necessary
background introduced in the previous sections,
here we will be finally able to state our {\it a priori} 
estimate for a solution $u\in W_{p,-1}^2(G)$ 
to the boundary value problem
\beqn\label{4.1}
&&\hskip -2truecm
\left\{\begin{array}{lll}
\l u(x) -{\cal A}(x;D_x)u(x) = f(x)\,, \qq x\in G\,,
\\[1mm]
{\cal B}(x;D_x)u(x)=g(x)\,, \qq x\in\partial G\,,
\end{array}
\right.
\eeqn
where $f\in L_{p,-1}(G)$, $g\in W_{p,-1}^{1-p^{-1}}(\partial G)$
and $\l\in\csp$.\\
However, to state the main result,
some basic assumptions on the domain $G$
and on the differential operators ${\cal A}(x;D_x)$ and ${\cal B}(x;D_x)$
are needed. We are going to list them.\\
Let $C_j$, $j=1,2,3$, be three positive constants
such that, denoted with $\phi$ the angle
$\arctan[(C_1)^{-1}]\in(0,\pi/2)$, they satisfy
$C_2>\sin\phi$ and $C_3>C_1$
and let $\eta:[0,C_2]\to\rsp$ be a function of class $C^2$
satisfying the following properties:
\begin{itemize}
\item[i)]\;$\eta(y)=C_1y$\,,\
{\rm if}\ $y\in [0,\sin\phi]$\,;
\item[ii)]\;$0<\eta(y)<C_3 y$\,,\
{\rm if}\ $y\in (\sin\phi,C_2)$\,;
\item[iii)]\;$\eta(C_2)=0$,\ $\eta'(C_2)=-\infty$\,.
\end{itemize}
Having such a $\eta$, for the rest of the paper
with $G$ we will denote the domain
\beqn\label{4.2}
&&\hskip -1truecm
G=\big\{x\in\rsp^n:\ |x'|< \eta(x_n)\,,\
x'=(x_1,\ldots,x_{n-1})\,,\ 0<x_n<C_2\big\}\,.
\eeqn
When $n\ge 2$, due to i), $G\cap B(0,1)$ coincides with the cone
$\{x\in\rsp^n:\! |x'|< C_1x_n,\ x_n>0\}$, whereas, when $n=1$,
we have $x'=0$ and
$G$ simply coincides with the interval $(0,C_2)$.
\begin{remark}\label{rem4.1}
\emph{To clarify the meaning of the constant $C_1$ and of the
choice of the interval $[0,\sin\phi]$ in the assumption i)
on $\eta$, observe first what happens when $n=2$.
In this case
$G\cap B(0,1)$ coincides with the cone $\{(x_1,x_2)\in\rsp^2:|x_1|<
C_1x_2, x_2>0\}$ and therefore,
using polar co-ordinates in $\rsp^2$ (set $\theta_0=\pi/2$ and $n=2$ in
formulae (\ref{3.4})), we deduce that for any $x\in G\cap B(0,1)$
the angle $\theta_1$ belongs to $(\phi,\pi-\phi)$.
Hence if $x\in \partial G\cap \ov{ B(0,1)}$ then $x_2\in [0,\sin\phi]$.
Generalizing to the case $n>2$,
by setting $\theta_0=\pi/2$ in formulae (\ref{3.4})
we deduce that for any $x\in G\cap B(0,1)$
all the angles $\theta_i$, $i=1,\ldots,n-1$, belong to the same interval
$(\phi,\pi-\phi)$.}
\end{remark}
Now, ${\cal A}(x;D_x)$ being defined by (\ref{1.1}), we assume
\beqn
\label{4.3}
&&\hskip -1truecm
a_{i,j}\in C^1(\ov G;\csp)\,,\q
a_{i,j}=a_{j,i}\,,\q
i,j=1,\dots,n\,,
\\[1mm]
\label{4.4}
&&\hskip -1truecm
{\rm Re}\sum_{i,j=1}^na_{i,j}(x)\xi_i\xi_j\ge C_0|\xi|^2\,,
\q\forall\,(x,\xi)\in (\ov G\backslash\{O\})\times \rsp^n\
{\rm and\ some}\ C_0>0\,.
\eeqn
As it is well-know, if $n\ge 2$ then
assumption (\ref{4.4}) implies the following:
\beqn
\label{4.5}
&&\hskip -1truecm
\left\{\begin{array}{lll}
for\ any\ x\in\ov G\ and\ any\ linearly\ independent\
vectors\ \xi,\zeta\in\rsp^n
\\[1mm]
the\ polynomial\;
\tau\to{\cal A}(x;\xi+\tau\zeta)
=\sum_{i,j=1}^na_{i,j}(x)(\xi_i+\tau\zeta_i)(\xi_j+\tau\zeta_j)
\\[1mm]
has\ a\ unique\ root\ with\ positive\ imaginary\ part.
\end{array}\right.
\eeqn
With ${\cal A}(x;D_x)$ we associate the boundary operator
\beqn\label{4.6}
&&\hskip -1,5truecm
{\cal B}(x;D_x)=\sum_{i=1}^nb_i(x)D_{x_i} + b_0(x)I\,,
\q x\in\partial G\,,
\eeqn
where, ${\cal V}$ being an open neighborhood of $\ov G$,
we have
\beqn\label{4.7}
&&\hskip -1,7truecm
b_j\in C^1({\cal V};\rsp)\,,\qq j=0,1,\ldots,n\,.
\eeqn
Moreover, if $n\ge 2$ we assume that 
the $b_j$'s satisfy also the
following two requirements:
\beqn\label{4.8}
&&\hskip -1,3truecm
\left\{\begin{array}{lll}
|\sum_{i=1}^{n-1}b_i(x'\cos\gamma,x_n)v_i(x)\cos\gamma
+b_n(x'\cos\gamma,x_n)v_n(x)+|v'(x)|^2\sin^2\g| \ge m>0\,,
\\[2mm]
for\ any\  \g\in[0,2\pi]\ and\ any\ x\in\partial G\backslash\{O\},
\\[2mm]
v(x)=(v'(x),v_n(x))
\ being\ the\ outer\ normal\ to\
\partial G\backslash\{O\}\ at\ x\,,
\end{array}
\right.
\eeqn
\beqn
&&\hskip -1,8truecm
\label{4.9}
\left\{\!\!\begin{array}{lll}
for\ any\ x\in\partial G\backslash\{O\}
\ and\ for\ any\ \xi, \zeta\in\rsp^n,
repectively\ tangent\ and\ normal
\\[1mm]
to\ \partial G\backslash\{O\}\ at\ x,
the\ polynomial\ {\cal B}(x;\xi+\tau\zeta)=
\sum_{i=1}^nb_i(x)(\xi_i+\tau\zeta_i)\ is\ not
\\[1mm]
divisible\ by\ (\tau-\tau^+(x,\xi,\zeta))\
without\ remainder,\ \tau^+(x,\xi,\zeta)\ being\ the\ unique
\\[1mm]
root\ with\ positive\ imaginary\ part\ of\ the\ 
polynomial\ {\cal A}(x;\xi+\tau\zeta)\ in\ (\ref{4.5})\,.
\end{array}
\right.
\eeqn
Observe that (\ref{4.8}) is well defined by virtue of (\ref{4.7})
since if $x\in\partial G\backslash\{O\}$ then $(x'\cos\gamma,x_n)$,
$\g\in[0,2\pi]$, belongs to $G$. Assumption (\ref{4.8})
can be considered as an improvement of the standard
assumption for the coefficients of ${\cal B}(x;D_x)$,
corresponding to $\g=0$ in (\ref{4.8}). In Section 5
we will exhibit a concrete class of functions $b_j$, $j=1,\ldots,n$, 
satisfying (\ref{4.8}).
\\
Instead, in the case $n=1$ we assume
\beqn\label{4.10}
b_1(x_1)\neq [\eta'(x_1)]^{-1}\eta(x_1),
\q\ {\rm for\ every}\, x_1\in (0,C_2]
\ {\rm such\ that}\ \eta'(x_1)\neq 0\,,
\eeqn
with the convention that when $x_1=C_2$
then (\ref{4.10}) should be understood as $b_1(C_2)\neq 0$.
\\
In accordance with Definition 1.5 on p. 113 in
\cite{LM}, assumption (\ref{4.9})
means that ${\cal B}(x;D_x)$ covers ${\cal A}(x;D_x)$
on $\partial G\backslash\{O\}$. 
We will need the following preliminary result.
\begin{lemma}\label{lem4.2}
Let $G$ be the domain defined by (\ref{4.2})
and let ${\cal A}(x;D_x)$ and ${\cal B}(x;D_x)$ be
the differential operators defined respectively
by (\ref{1.1}) and (\ref{4.6}) and assume that the coefficients of
${\cal A}(x;D_x)$ satisfy (\ref{4.3}), (\ref{4.4}) whereas the
coefficients of ${\cal B}(x;D_x)$ satisfy (\ref{4.7})--(\ref{4.10}).
Denote with $(x_0,x)$, $x=(x',x_n)$,
the points of $\rsp^{n+1}$
and with $\wtil G$ the domain
\beqn\label{4.11}
&&\hskip -1truecm
\wtil G=
\big\{(x_0,x)\in\rsp^{n+1}:
|(x_0,x')|< \eta(x_n)\,,\ 0< x_n<C_2
\big\}.
\eeqn
Then ${\cal B}((x_0,x);D_x, D_{x_0})={\cal B}(x;D_x)+x_0D_{x_0}$ covers
${\cal A}_\psi(x;D_x,D_{x_0})={\cal A}(x;D_x)+e^{i\psi}D_{x_0}^2$,
$\psi\in[-\pi/2,\pi/2]$, on $\partial \wtil G\backslash\{O\}$.
\end{lemma}
\pn
The proof Lemma \ref{lem4.2} will be given in Section 5. 
Here we make only two easy remarks.
\begin{remark}\label{rem4.3}
\emph{Observe that if $n=2$
then $\wtil G$ is a 3--dimensional domain
generated by a rotation of $G$ around the $x_2$--axis.
In this sense, when $n\ge2$, the $(n+1)$--dimensional domain
$\wtil G$ can always be viewed as a rotation of the
$n$--dimensional domain $G$ around the $x_n$--axis.
In addition, in the critical case $n=1$,
definition (\ref{4.11}) ensures that
$\wtil G$ is 2--dimensional domain, symmetric with respect to the $x_1$--axis
and coinciding with a cone near the origin. 
This will be important in
the following, since we will use Theorem \ref{thm3.2} in dimension $n+1$
and therefore we will need to consider $(n+1)$--dimensional
domain having the properties for which that theorem is true.}
\end{remark}
\begin{remark}\label{rem4.4}
\emph{Under our assumptions on ${\cal A}(x;D_x)$
the operator ${\cal A}_{\psi}(x;D_x,D_{x_0})$, 
$\psi\in[-\pi/2,\pi/2]$, does not necessarily satisfy (\ref{4.4}), 
but only the weaker ellipticity condition
\beqn
\label{4.12}
\Big|\sum_{i,j=1}^{n}a_{i,j}(x)\xi_i\xi_j
+e^{i\psi}\xi_0^2\Big|
\ge2^{-1/2}\min\{C_0,1\}|\xi|^2\,\q\;
\forall\,(x,\xi)\in(\ov G\backslash\{O\})\times \rsp^{n+1}.
\eeqn
However, if $n\ge 2$ then $n+1\ge 3$ and it
is well-known that in this case (\ref{4.12})
implies (\ref{4.5}).
If $n=1$, by computing explicitly the roots of the
polynomial ${\cal A}_\psi(x_1;\xi+\tau\zeta)$ for the operator
${\cal A}_\psi(x_1;D_{x_1},D_{x_0})$,
it can be checked that (\ref{4.5}) is satisfied, too.
Indeed, by virtue of (\ref{4.4}), when $n=1$ we may for simplicity assume
$a_{1,1}$ to be real and positive
and hence, given two
linearly independent vectors $\xi=(\xi_0,\xi_1),
\zeta=(\zeta_0,\zeta_1)\in\rsp^2$, it follows that the
polynomial ${\cal A}_\psi(x_1;\xi+\tau\zeta)=
a_{1,1}(x_1)(\xi_1+\tau\zeta_1)^2+e^{i\psi}(\xi_0+\tau\zeta_0)^2$
has the roots:
\beqn\label{4.13}
&&\hskip -2,5truecm
\tau_{\pm}(x_1,\xi,\zeta)=
[\zeta_0^4+(a_{1,1}(x_1))^2\zeta_1^4
+2a_{1,1}(x_1)\zeta_0^2\zeta_1^2\cos\psi]^{-1}
\no\\[1mm]
&&\hskip 0truecm
\times\Big\{-\xi_0\zeta_0^3-(a_{1,1}(x_1))^2\xi_1\zeta_1^3
-a_{1,1}(x_1)\zeta_0\zeta_1e^{-i\psi}
[\xi_0\zeta_1e^{2i\psi}+\xi_1\zeta_0]
\no\\
&&\hskip 0,6truecm
\pm\,i\,[a_{1,1}(x_1)\xi_0\zeta_1-\xi_1\zeta_0]e^{-i\psi/2}
[a_{1,1}(x_1)\zeta_1^2e^{i\psi}+\zeta_0^2]\Big\}.
\eeqn}
\end{remark}
Now, triplet $\{{\cal A}_{\psi}(x;D_x,D_{x_0}),
{\cal B}((x_0,x);D_x,D_{x_0}); \wtil G\}$
being defined as in Lemma \ref{lem4.2},
with problem (\ref{4.1}) we associate the following
boundary value problem, where ${\cal F}\in L_{p,-1}(\wtil G)$ and
${\cal G}\in W_{p,-1}^{1-p^{-1}}(\partial\wtil G)$
\beqn\label{4.14}
&&\hskip -2truecm
\left\{\begin{array}{lll}
{\cal A}_{\psi}(x;D_x,D_{x_0})v(x_0,x) = {\cal F}(x_0,x)\,,
\qq (x_0,x)\in \wtil G\,,
\\[1mm]
{\cal B}((x_0,x);D_x,D_{x_0})v(x_0,x)={\cal G}(x_0,x)\,, \qq
(x_0,x)\in\partial\wtil G\,,
\end{array}
\right.
\eeqn
As shown in Remark \ref{rem3.1} when the coefficients $a_{i,j}$,
$i,j=1,\ldots,n$, of ${\cal A}(x;D_x)$
satisfy (\ref{4.3}) then ${\cal A}_{\psi}(x;D_x,D_{x_0})$
belongs to the class ${\frak C}(2,0)$. Moreover, when the
coefficients $b_h$, $h=0,1,\ldots,n$, of
${\cal B}(x;D_x)$ satisfy (\ref{4.7}) then,
using formulae (\ref{3.4}), (\ref{3.5}) and (\ref{3.15}),
it is easy to see that ${\cal B}((x_0,x);D_x,D_{x_0})$
belongs to ${\frak C}(1.1)$.
In addition formula (\ref{4.12}) shows that
${\cal A}_{\psi}(x;D_x,D_{x_0})$ is uniformly elliptic in
$\wtil G\backslash\{O\}$, whereas Lemma \ref{lem4.2}
establishes that the boundary conditions
${\cal B}((x_0,x);D_x, D_{x_0})$ covers
${\cal A}_\psi(x;D_x,D_{x_0})$ on $\partial \wtil G\backslash\{O\}$.
\\
Hence, denoted with $\Om$ the intersection $\wtil G\cap \partial B(0,1)$,
 problem (\ref{4.14}) generates a model problem
in the cone $(0,+\infty)\times \Om$.
Indeed, $\theta_i$, $i=0,\ldots,n-1$,
being defined by (\ref{3.4}), with the pair
$\{{\cal A}_\psi(x;D_x,D_{x_0}),{\cal B}((x_0,x);D_x, D_{x_0})\}$
we associate the operator
\beqn\label{4.15}
&&\hskip -1truecm
{\cal U}(0,\om,z, D_{\om})=\{A_{\psi}(0,\om,z-2i,D_{\om}),
B(0,\om,z-2i,D_{\om})\}\,,
\eeqn
where $\om=(\theta_0,\ldots,\theta_{n-1})$,
$D_{\om}=(\theta_0,\ldots,\theta_{n-1})$,
and the operators $A_{\psi}$ and $B$ are determined
from ${\cal A}_\psi(x;D_x,D_{x_0})$ and ${\cal B}((x_0,x);D_x, D_{x_0})$
by means of (\ref{3.1}) replacing the coefficients $p_{h,\a}(r,\om)$
with $p_{h,\a}(0,\om)$.
Since problem (\ref{4.14}) is uniformly elliptic in $\wtil G
\backslash\{O\}$ and the boundary condition covers
${\cal A}_\psi(x;D_x,D_{x_0})$ on $\partial \wtil G\backslash\{O\}$
then the operator ${\cal U}(0,\om,z, D_{\om})$
is elliptic with complex parameter in the sense of \cite{AV}.
Therefore, with the choice of the parameter as in (\ref{3.22})
and $G$ replaced by $\wtil G$,
all the assumptions on problem (\ref{3.16})
which are necessary in order to state Theorem \ref{thm3.2}
are satisfied even for problem (\ref{4.14}).\\
We now give the following definition,
arising from the necessity to use Theorem \ref{thm3.2}
in dimension $n+1$, with $\b=-1$ and $l=0$.
\begin{definition}\label{def4.5}
The boundary value
problem (\ref{4.14}) will be said
{\rm regular} if the line ${\rm Im}\, z=-1+(n+1)p^{-1}$
contains no eigenvalue of the operator ${\cal U}(0,\om,z, D_{\om})$
defined by (\ref{4.15}).
In this case the triplet $\{{\cal A}(x;D_x),{\cal B}(x;D_x); G\}$
will be said the restriction to the $x$ variable of the
regular boundary value problem (\ref{4.14}) in
the domain $\wtil G$ related to $G$ by (\ref{4.11}).
\end{definition}
 
Definition \ref{def4.5} can be considered as the
equivalent, in the setting of weighted Sobolev function spaces,
of Definition 6.2 in \cite{AN}.
In this sense our results are
in accordance with those proven in \cite{AN}  for
the subclass of problems consisting
in the restriction to the $x$ variable of regular
elliptic problems in one more variable.
Two simple examples for the Definition \ref{def4.5}
are those given on p. 45 in \cite{MP1}
and p. 86 in \cite{MP2} and related to the
homogeneous Dirichlet boundary conditions.\newpage\pn
We can finally state our main result.
\begin{theorem}\label{thm4.6}
Let $p>n$, $p\neq n+1$, and let
the triplet $\{{\cal A}(x;D_x),{\cal B}(x;D_x); G\}$
be the restriction to the $x$ variable of the regular boundary
value problem (\ref{4.14}) in the sense of Definition \ref{def4.5}.
Then, $g_0\in W_{p,-1}^1(G)$ being any extension
to $G$ of ${\cal B}(x;D_x)u$,
there exists $\om>0$ such that if ${\rm Re}\l\ge\om$
for every $u\in W_{p,-1}^2(G)$ the following estimate hold
\beqn\label{4.16}
&&\hskip -1truecm
|\l|\|u\|_{L_{p,-1}(G)}
+|\l|^{1/2}\|Du\|_{L_{p,-1}(G)}
+\|D^2u\|_{L_{p,-1}(G)}
\no\\[2mm]
&&\hskip -1,5truecm
\le M\big\{
\|\big(\l I-{\cal A}(x;D_x)\big) u\|_{L_{p,-1}(G)}
+\,(1+|\l|^{1/2})\|g_0\|_{L_{p,-1}(G)}+\|Dg_0\|_{L_{p,-1}(G)}\big\}.
\eeqn
The positive constant $M$ in (\ref{4.16}) depends only on $p$, $n$, the
$C^1(G)$-norm of the coefficients of ${\cal A}(x;D_x)$
and of ${\cal B}(x;D_x)$ and the constants $C_j$, $j=2,3$, intervening
in the properties i)--iii) for the function $\eta$ which
describes the boundary $\partial G$ of $G$.
\end{theorem}
\pn As announced in the Introduction, the proof of Theorem \ref{4.6}
will be given in Section 7. Here, instead, we want to show
some easy consequence of estimate (\ref{4.16}). We set
\beqn\label{4.17}
&&\hskip -2truecm
\left\{
\begin{array}{lll}
{\cal D}(A)=\{u\in W_{p,-1}^2(G):\ {\cal B}(x;D_x)u=0\ \;
{\rm in}\;\partial G\}\,,
\\[1mm]
Au={\cal A}(x;D_x)u\,,\qq\forall\,u\in{\cal D}(A)\,.
\end{array}
\right.
\eeqn
$A$ is said the realization of ${\cal A}(x;D_x)$ in  $L_{p,-1}(G)$
with homogeneous boundary condition.
\begin{corollary}\label{cor4.7}
Let assumptions of Theorem \ref{thm4.6} be fulfilled
and let the pair $(A,{\cal D}(A))$ be defined by (\ref{4.17}).
There exists $\om>0$ such that if ${\rm Re}\l\ge\om$
then the operator $\l I-A$ is closed and injective in  $L_{p,-1}(G)$. 
As a consequence,
$A$ is closed in  $L_{p,-1}(G)$.
\end{corollary}
\begin{proof}
By taking as $g_0$ the null function, the injectivity 
of $\l I-A$, ${\rm Re}\l\ge\om$,
trivially follows from (\ref{4.16}). Now, let
$\{u_n\}_{n\in\nsp}\subset{\cal D}(A)$ 
such that $u_n\to u$ in $L_{p,-1}(G)$ and
$(\l I-A)u_n\to v$ in $L_{p,-1}(G)$. If we set $g_0=0$, from (\ref{4.16}) 
it clearly follows that $\{u_n\}_{n\in\nsp}$ is a Cauchy
sequence in $W_{p,-1}^2(G)$ and hence $u\in W_{p,-1}^2(G)$.
Moreover, due to assumptions (\ref{4.3}), (\ref{4.7})
on the coefficients of ${\cal A}(x;D_x)$ and ${\cal B}(x;D_x)$, respectively,
it is easy to deduce $(\l I-A)u_n\to(\l I-A)u $ in $L_{p,-1}(G)$
and $0={\cal B}(x;D_x)u_n\to{\cal B}(x;D_x)u$ 
in $W_{p,-1}^{1-p^{-1}}(\partial G)$. Therefore $u\in {\cal D}(A)$
and $(\l I-A)u=v$, i.e. $\l I-A$ is closed in $L_{p,-1}(G)$.\\
The last assertion follows from $A=\l I -(\l I-A)$.
\end{proof}
As a further corollary of Theorem \ref{thm4.6} we show that, if
a solution  $u\in W_{p,-1}^2(G)$ of problem (\ref{4.1}) exists, then
the operator $A$ defined via (\ref{4.17}) is sectorial, 
in accordance with the Definition 2.0.1 in \cite{LU} 
which we report for reader's convenience.
\begin{definition}\label{def4.8}
Let $X$ be a complex Banach space, with norm $\|\cdot\|$. A linear operator
$B:{\cal D}(B)\subset X\to X$ is said to be {\it sectorial} if there
are constants $\om_1\in\rsp$, $\vartheta\in(\pi/2,\pi)$, $C>0$ such that,
denoted with $\rho(B)$ the resolvent set of $B$, the following hold:
\\[1mm]
(i)\ \ $\rho(B)\supset S_{\vartheta,\om_1}=
\{z\in\csp:\ z\neq\om_1,|\arg(z-\om_1)|<\vartheta\},$
\\[2mm]
(ii)\ $\|(zI-B)^{-1}\|_{{\cal L}(X)}\le C|z-\om_1|^{-1},\q\forall \,z\in 
S_{\vartheta,\om_1}.$
\end{definition}
\pn We recall also the following sufficient condition for an operator
to be sectorial and for the proof of which we refer to
\cite[Proposition 2.1.11]{LU}.
\begin{proposition}\label{prop4.9}
Let $\om_1\in\rsp$ and let $B:{\cal D}(B)\subset X\to X$ 
be a linear operator such that
$\rho(B)\supset S_{\om_1}=\{z\in\csp:\ {\rm Re}z\ge\om_1\}$ and
$\|(z I -B)^{-1}\|_{{\cal L}(X)}\le C|z|^{-1}$ 
for any $z\in S_{\om_1}$ and some $C>0$. Then $B$ is sectorial.
\end{proposition}
\pn Consequently, we have the following corollary.
\begin{corollary}\label{cor4.10}
Let assumptions of Theorem \ref{thm4.6} be fulfilled and let us suppose
that for any pair $(f,g)\in L_{p,-1}(G)\times W_{p,-1}^{1-p^{-1}}(\partial G)$
there exists a solution $u\in W_{p,-1}^2(G)$ to problem (\ref{4.1}).
Then the operator $A$ defined by (\ref{4.17}) is sectorial.
\end{corollary}
\begin{proof}
If for any pair $(f,g)\in L_{p,-1}(G)\times W_{p,-1}^{1-p^{-1}}(\partial G)$
a solution $u\in W_{p,-1}^2(G)$ of problem (\ref{4.1}) exists,
then, when ${\rm Re}\l\ge\om$, the solution is unique
by virtue of estimate (\ref{4.16}). Moreover, from (\ref{4.16}) with
$g_0$ equal to zero, we deduce
$\|(\l I -A)^{-1}\|_{{\cal L}(L_{p,-1}(G))}\le M|\l|^{-1}$
for any $\l\in\csp$ such that ${\rm Re}\l\ge\om$.
Hence, the assertion follows from Proposition \ref{prop4.9}.
\end{proof}
Since it is well-known that sectorial operators 
generate analytic semigroups, we have the following
further corollary.
\begin{corollary}\label{cor4.11}
Under the hypotheses of Corollary \ref{cor4.10} 
the realization $A$ of ${\cal A}(x;D_x)$ 
in $L_{p,-1}(G)$ with homogeneous boundary condition 
generates an analytic semigroup of linear bounded operators 
$\{{\cal T}(t)\}_{t\ge 0}\subset{\cal L}(L_{p,-1}(G))$.
\end{corollary}
\section{Proof of Lemma \ref{lem4.2}}
\setcounter{equation}{0}
First, accordingly to Remark {\ref{rem4.3},
we observe that if $n\ge 2$ and $\wtil G$
is related to $G$ by (\ref{4.11}) then
a very special characterization
of the points in $\partial\wtil G\backslash\{O\}$
in terms of those in $\partial G\backslash\{O\}$ can be given.
Indeed, when $n\ge 2$ and
$(\wtil x_0,\wtil x)\in\partial\wtil G\backslash\{O\}$
(i.e. $|(\wtil x_0,\wtil x')|=\eta(\wtil x_n)$) 
we set $\a=|(\wtil x_0,\wtil x')|^{-1}|\wtil x'|$,
$\b=|(\wtil x_0,\wtil x')|^{-1}\wtil x_0$. 
Since $\a^2+\b^2=1$ there exists $\varphi\in [0,2\pi]$ 
such that $\a=|\cos\varphi|$ and $\b=\sin\varphi$.
Let us set $x_i=\wtil x_i/\cos\varphi$, $i=1,\ldots,n-1$, $x_n=\wtil x_n$.
If $\cos\varphi=0$, i.e. when $(\wtil x_0,\wtil x)$ is of the form
$(\wtil x_0,0,\ldots,0,\wtil x_n)$, we set $x_i=0$, $i=1,\ldots,n-1$,
$x_n=C_2$. In this way we have defined a point
$x\in\partial G\backslash\{O\}$. In fact,
$|x'|=|\wtil x'|/|\cos\varphi|=|\wtil x'|/\a=|(\wtil x_0,\wtil x')|$
if $\a\neq 0$ and $x=(0,\ldots,0,C_2)$ if $\a=0$.
Summing up, if $n\ge 2$, given 
$(\wtil x_0,\wtil x)\in\partial\wtil G\backslash\{O\}$
there exists an angle $\varphi\in[0,2\pi]$ such that
\beqn\label{5.1}
&&\hskip -1truecm
(\wtil x_0,\wtil x',\wtil x_n)=
(|x'|\sin\varphi,x'\cos\varphi, x_n)\,,\q
x=(x',x_n)\in\partial G\backslash\{O\}\,,
\eeqn
where boundary points of the form $(\wtil x_0,0,\ldots,0,\wtil x_n)$
correspond to the choice $\varphi=\pi/2$ if $\wtil x_0>0$ and
$\varphi=3\pi/2$ if $\wtil x_0<0$.\\ \vskip -0,2truecm \pn
{\it Proof of Lemma \ref{lem4.2}.}
We consider the following two distinct cases:
{\it i)} $n\ge 2$, {\it ii)} $n=1$.
\\
{\it i)} $n\ge 2$.
Let $\Phi:\rsp^{n+1}\to \rsp$ to be the function defined by
\beqn\label{5.2}
&&\hskip -2truecm
\Phi(y_0,y)={2}^{-1}\{|(y_0,y')|^2-[\eta(y_n)]^2\}\,,\q y_n>0\,.
\eeqn
It follows
$\partial\wtil G\backslash\{O\}=
\{(y_0,y)\in\rsp^{n+1}:\Phi(y_0,y)=0,\,y_n>0\}$
and hence, since the only point of
$\partial \wtil G\backslash\{O\}$ with $x_n=C_2$
is the point $(0,\ldots,0,C_2)$ with normal
$(0,\ldots,0,1)$, the normal $\wtil\zeta$ to
$\partial\wtil G\backslash\{O\}$ at $(\wtil x_0,\wtil x)$
is given by
\beqn\label{5.3}
&&\hskip -1truecm
\wtil \zeta=
\left\{\!
\begin{array}{lll}
(\wtil x_0,\wtil x',-\eta(\wtil x_n)\eta'(\wtil x_n))\,,
\q\; {\rm if}\;\;x_n\in(0,C_2)\,,
\\[1mm]
(0,\ldots,0,1)\,,
\q\;{\rm if}\;\;x_n=C_2\,.
\end{array}
\right.
\eeqn
From (\ref{5.2}) it follows also that the normal $v(x)$ to
$\partial G\backslash\{O\}$ at $x$ is the vector
\beqn\label{5.4}
&&\hskip -1,5truecm
v(x)=\left\{\!
\begin{array}{lll}
(x',-\eta(x_n)\eta'(x_n))\,,
\q\; {\rm if}\;\;x_n\in(0,C_2)\,,
\\[1mm]
(0,\ldots,0,1)\,,
\q\;{\rm if}\;\;x_n=C_2\,,
\end{array}
\right.
\eeqn
and, since $v'(x)=x'$, we see
that (\ref{5.3})
can be rewritten in the more compact way
\beqn\label{5.5}
&&\hskip -2truecm
\wtil \zeta=
(|v'(x)|\sin\varphi,v'(x)\cos\varphi,v_n(x))\,,
\eeqn
where $v(x)$ defined by (\ref{5.4})
is the normal at the point $x$ such that (\ref{5.1}) holds.\\
It remains to characterize the tangent vectors $\wtil\xi$ to
$\partial\wtil G\backslash\{O\}$.
Taking advantage from (\ref{5.4}), (\ref{5.5})
it is not too difficult to show
that any vector $\wtil\xi$ tangent
to $\partial \wtil G\backslash\{O\}$
at $(\wtil x_0,\wtil x)$ has one
of the following three representations
\beqn\label{5.6}
&&\hskip -1truecm
\left\{\!\!
\begin{array}{lll}
\big(a\cos\varphi+c|v'(x)|\sin\varphi,
by'+cv'(x)\cos\varphi-a\,\ds\frac{v'(x)}{|v'(x)|}\sin\varphi ,
c\,\ds\frac{\eta(x_n)}{\eta'(x_n)}\,\big),\q{\rm if}\;\,\eta'(x_n)\neq 0,
\\[3mm]
\big(a\cos\varphi,by'-a\,\ds\frac{v'(x)}{|v'(x)|}\sin\varphi,c\big),
\q{\rm if}\;\,\eta'(x_n) =0\;\; {\rm and}\ \;\varphi\notin\{\pi/2,3\pi/2\}
\\[3mm]
(0,y),\; y\in \rsp^n,\; |y|\neq 0,
\q{\rm if}\;\,\eta'(x_n) =0\;\; {\rm and}\ \;\varphi\in\{\pi/2,3\pi/2\}
\end{array}
\right.
\eeqn
where $y'\in\rsp^{n-1}$, $|y'|\neq 0$, satisfies
$y'\cdot x'=0$ if $x_n\neq C_2$ or $\varphi\notin\{\pi/2,3\pi/2\}$,
$a,b,c$ are real numbers not all equal to zero and
$x\in\partial G\backslash\{O\}$ is the point in (\ref{5.1}).\\
Now, let $(\wtil x_0,\wtil x)\in\partial \wtil G\backslash\{O\}$
with $\wtil x_0\neq 0$
and assume
condition (\ref{4.9}) is violated for the operator
${\cal B}((\wtil x_0,\wtil x);D_x,D_{x_0})$.
Hence, denoted with
$\wtil\tau^+((\wtil x_0,\wtil x),\wtil\xi,\wtil\zeta)$ the
unique root with positive imaginary part of the polynomial
${\cal A}_{\psi}(\wtil x,\wtil\xi+\tau\wtil\zeta)$
we have, for any $\tau\in\csp$
\beqn\label{5.7}
\hskip -1truecm
&&\!\!
\sum_{i=1}^{n}b_i(\wtil x)(\wtil\xi_i+\tau\wtil\zeta_i)+
\wtil x_0(\wtil\xi_0+\tau\wtil\zeta_0)=
\chi((\wtil x_0,\wtil x),\wtil\xi,\wtil\zeta)
[\tau-\wtil\tau^+((\wtil x_0,\wtil x),\wtil\xi,\wtil\zeta)]\,.
\eeqn
From (\ref{5.6}) we deduce that there are
three different situations to take into examination.\\
1) {\it Case $\eta'(x_n)\neq 0$}.
In this case from (\ref{5.1}),
(\ref{5.5})--(\ref{5.7}) and assumption (\ref{4.8})
we easily deduce
\beqn\label{5.8}
&&\hskip -0,5truecm
\sum_{i=1}^{n-1}b_i(x'\cos\varphi,x_n)v_i(x)\cos\varphi
+b_n(x'\cos\varphi,x_n)v_n(x)+|v'(x)|^2\sin^2\g\no
\\
&&\hskip -0,8truecm
=\,\chi((\wtil x_0,\wtil x),\wtil\xi,\wtil\zeta)\neq 0\,,
\\[2mm]
\label{5.9}
&&\hskip -0,5truecm
\Big\{\sum_{i=1}^{n-1}b_i(x'\cos\varphi,x_n)[by'_i+
cv_i(x)\cos\varphi-a|v'(x)|^{-1}v_i(x)\sin\varphi]
\no\\
&&\hskip -0,2truecm
+\,c\,b_n(x'\cos\varphi,x_n)[\eta'(x_n)]^{-1}\eta(x_n)
+|v'(x)|\big(a\cos\varphi+c|v'(x)|\sin\varphi\big)\sin\varphi\Big\}
\no\\[1mm]
&&\hskip -0,8truecm
=\,-\chi((\wtil x_0,\wtil x),\wtil\xi,\wtil\zeta)\,
\wtil\tau^+((\wtil x_0,\wtil x),\wtil\xi,\wtil\zeta)
\eeqn
From (\ref{5.8}), (\ref{5.9}) and the fact that
$b_j$, $j=1,\ldots,n$, assume only real values (cf. (\ref{4.7}))
we get a contradiction
since ${\rm Im}\,\tau^+((\wtil x_0,\wtil x),\wtil\xi,\wtil\zeta)>0$.\\
2) {\it Case $\eta'(x_n)=0$, $\varphi\notin\{\pi/2,3\pi/2\}$}.
From (\ref{5.4}) we see that in this case we have $v_n(x)=0$
and hence from (\ref{5.8}) the contradiction follows as in the
case before using assumption (\ref{4.8}) and changing
the left-hand side of (\ref{5.9}) in accordance with (\ref{5.6}).\\
3) {\it Case $\eta'(x_n)=0$, $\varphi\in\{\pi/2,3\pi/2\}$}.
From (\ref{5.1}), (\ref{5.3}) and (\ref{5.6})
we obtain that $\wtil\zeta$ and $\wtil\xi$
are respectively of the form $(\wtil\zeta_0,0,\ldots,0)$
and $(0,y)$ with  $\wtil\zeta_0\neq 0$ and $y\in\rsp^n$, $|y|\neq 0$.
However, due to the assumptions on $\eta$ it follows $x_n\in (0,C_2)$ and
hence, since $(\wtil x_0,0,\ldots,0,x_n)\in\partial\wtil G\backslash\{O\}$,
we have $\wtil x_0=\eta(x_n)\neq 0$. Therefore, from (\ref{5.7}) we get
\beqn\label{5.10}
&&\hskip -1,7truecm
\wtil x_0\wtil\zeta_0
=\chi((\wtil x_0,0,\ldots,0,x_n),\wtil\xi,\wtil\zeta)\neq 0\,,
\\
\label{5.11}
&&\hskip -1,7truecm
\sum_{i=1}^{n}b_i(0,\ldots,0,x_n)y_i
=-\chi((\wtil x_0,0,\ldots,0,x_n),\wtil\xi,\wtil\zeta)\,
\wtil\tau^+((\wtil x_0,0,\ldots,0,x_n),\wtil\xi,\wtil\zeta)\,,
\eeqn
and again the contradiction follows from the fact that
the $b_i$'s assume only real values whereas
${\rm Im}\,\tau^+((\wtil x_0,\wtil x),\wtil\xi,\wtil\zeta)>0$.\\
Contradictions we get in 1)--3) mean that the assumption 
that condition (\ref{4.9}) was violated for
${\cal B}((\wtil x_0,\wtil x);D_x,D_{x_0})$ was wrong
and so, if $n\ge 2$,
the proof is complete.
\\
{\it ii)} $n=1$. In this case, since $x'=0$, 
no relationship of type (\ref{5.1})
is possible and we can not reason
as before. However, since the points of
$\partial\wtil G\backslash\{O\}$
are the points $(\wtil x_0,\wtil x_1)=(\pm \eta(x_1),x_1)$,
from (\ref{5.2}) with $n=1$
we deduce that the normal $\wtil\zeta$ and the tangent $\wtil\xi$ to
$\partial\wtil G\backslash\{O\}$ at $(\pm \eta(x_1),x_1)$
have the following form
\beqn\label{5.12}
&&\hskip -1,5truecm
\wtil\zeta=
\left\{\!\!
\begin{array}{lll}
(\pm \eta(x_1),-\eta(x_1)\eta'(x_1))\,,\q{\rm if}\;\;x_1\in(0,C_2)\,,
\\[2mm]
(0,1)\,,\q{\rm if}\;\;x_1=C_2\,,
\end{array}
\right.
\\[2mm]
\label{5.13}
&&\hskip -1,5truecm
\wtil\xi=
\left\{\!\!
\begin{array}{lll}
(\pm \eta(x_1),[\eta'(x_1)]^{-1}\eta(x_1))\,,
\q{\rm if}\;\;\eta'(x_1)\neq 0\,,
\;{\rm and}\;\ x_1\in(0,C_2)\,,
\\[2mm]
(0,1)\,,
\q{\rm if}\;\;\eta'(x_1)=0\,,
\\[2mm]
(1,0)\,,
\q{\rm if}\;\;x_1=C_2\,.
\end{array}
\right.
\eeqn
Now, assume that (\ref{4.9}) does not hold for the operator
${\cal B}((\wtil x_0,\wtil x);D_x,D_{x_0})$. Hence, denoted with
$\wtil\tau^+((\wtil x_0,\wtil x_1),\wtil\xi,\wtil\zeta)$ the
unique root with positive imaginary part of the polynomial
${\cal A}_{\psi}(\wtil x,\wtil\xi+\tau\wtil\zeta)$ (cf. (\ref{4.13}))
we have that (\ref{5.7}) reduces to
\beqn\label{5.14}
&&\hskip -1truecm
b_1(x_1)(\wtil\xi_1+\tau\wtil\zeta_1)+\wtil x_0(\wtil\xi_0+\tau\wtil\zeta_0)=
\chi((\wtil x_0,\wtil x),\wtil\xi,\wtil\zeta)
[\tau-\wtil\tau^+((\wtil x_0,\wtil x),\wtil\xi,\wtil\zeta)]\,.
\eeqn
From (\ref{5.12}), (\ref{5.13}) we deduce that only
two situations have to be examined.\\
1) {\it Case} $x_1\in (0,C_2)$.
If $\eta'(x_1)\neq 0$, from (\ref{5.12})--(\ref{5.14})
and assumption (\ref{4.10}) we find
\beqn\label{5.15}
&&\hskip -1truecm
-b_1(x_1)\eta(x_1)\eta'(x_1)+[\eta(x_1)]^2=
\chi((\wtil x_0,\wtil x),\wtil\xi,\wtil\zeta) \neq 0
\\[2mm]
\label{5.16}
&&\hskip -1truecm
b_1(x_1)[\eta'(x_1)]^{-1}\eta(x_1)+[\eta(x_1)]^2
=-\chi((\wtil x_0,\wtil x),\wtil\xi,\wtil\zeta)
\wtil\tau^+((\wtil x_0,\wtil x),\wtil\xi,\wtil\zeta)
\eeqn
which is a contradiction since on the left-hand side of (\ref{5.16})
we have a real value whereas on the right-hand side we have
a complex number with positive imaginary part.
It is easy to observe that if $\eta'(x_1)= 0$
we still get a contradiction. Indeed,
due to the fact that $x_1\in (0,C_2)$,
on the left-hand side of (\ref{5.15})
we have $[\eta(x_1)]^2\neq 0$  whereas (cf. (\ref{5.13}))
on the left-hand side of (\ref{5.16})
we have only the real value $b_1(x_1)$.\\
2) {\it Case} $x_1= C_2$.
Since $\eta(C_2)=0$, from (\ref{5.12})--(\ref{5.14}) we find
\beqn\no
&&\hskip -1truecm
b_1(C_2)\tau=
\chi((0,C_2),\wtil\xi,\wtil\zeta)
[\tau-\wtil\tau^+((0,C_2),\wtil\xi,\wtil\zeta)]\,,
\eeqn
which is a contradiction due to assumption (\ref{4.10}).
Hence, also in the case $n=1$ we are done and the proof of
Lemma \ref{lem4.2} is complete.
\begin{remark}\label{rem5.1}
\emph{
With the help of (\ref{5.4})
we present here a class of coefficients $b_j$, $j=1,\ldots,n$,
which satisfy assumption (\ref{4.8}).
To this purpose, for any $x\in\ov G$ let set
\beqn\label{5.17}
&&\hskip -2truecm
b_j(x)=x_j\,,\q j=1,\ldots,n-1,\qq
b_n(x)\neq [\eta'(x_n)]^{-1}\eta(x_n)\,,\q{\rm if}\;\;\eta'(x_n)\neq 0.
\eeqn
Since from (\ref{5.4}) it follows
$v'(x)=x'$ for any $x\in \partial G\backslash\{O\}$,
with the coefficients
defined by (\ref{5.17}) and  using $|x'|=\eta(x_n)$
we see that (\ref{4.8}) is equivalent to require
\beqn
&&\hskip -1,5truecm
[\eta(x_n)]^2-b_n(x'\cos\g,x_n)\eta(x_n)\eta'(x_n)\neq 0\,\qq
{\rm if}\;\; x_n\in (0,C_2)\,,\no
\\[1mm]
&&\hskip -1,5truecm
-b_n(0,\ldots,0,C_2)\neq 0\qq {\rm if}\;\;x_n\in (0,C_2)\,.
\no
\eeqn
Therefore, with the convention that the assumption on $b_n$ in (\ref{5.17})
should be intended as $b_n(x)\neq 0$ if $x_n=C_2$ (i.e. when
$\eta'(x_n)=-\infty$), the previous two
inequality are both satisfied even in the case $\eta'(x_n)=0$
since in this case we have $x_n\in (0,C_2)$ and hence
$\eta(x_n)\neq0$. Observe also that in the case
$n=1$ then (\ref{5.17}) corresponds to (\ref{4.10}).}
\end{remark}
\begin{remark}\label{rem5.2}
\emph{
From (\ref{5.10}) and (\ref{5.11})
we see that in the case $n>1$
Lemma \ref{lem4.2} fails if
instead of ${\cal B}((x_0,x);D_x,D_{x_0})$
we consider only the operator ${\cal B}(x;D_x)$.
In the case $n=2$, using complex valued coefficients $b_j$, $j=1,2$,
there could be still the possibility to conclude the proof considering
only ${\cal B}(x;D_x)$, but surely no  if $n\ge 3$.}
\end{remark}
\section{The cut-off function}
\setcounter{equation}{0}
The procedure we will perform in Section 7 to prove estimate (\ref{4.16})
for a function $u\in W_{p,-1}^2(G)$, $G$ being defined by (\ref{4.2}),
requires the implementing of a function $v$
depending on $n+1$ variables
and having the following form
\beqn\label{6.1}
&&\hskip -1truecm
v(x_0,x)=\varkappa(x_0,x)e^{i\rho x_0}u(x)\,,\qq(x_0,x)\in\wtil G\,,
\eeqn
where $\rho>0$, $\wtil G$ is related to $G$ by (\ref{4.11})
and $\varkappa$ is an infinitely differentiable function
having compact support on $\wtil G$.\\
Functions of type (\ref{6.1}), with the aim
of proving estimates for the function $u$, are used, for instance,
in \cite{A}, \cite{AN} and \cite{LU}.
However, in that papers the domain $\wtil G$
always consists in the infinite ``cylinder'' $\G=(-\infty,+\infty)\times G$
and this, as remarked in the Introduction,
allows the authors to use cut-off functions $\varkappa$
depending only on $x_0$. \\
In our case the situation is really different, since when $\wtil G$
is related to $G$ by (\ref{4.11}) then $(x_0,x)$
do not belongs to $\wtil G$ if $|x_0|>\{[\eta(x_n)]^2-|x'|^2\}^{1/2}$.
By recalling definition of $\phi$ before the definition (\ref{4.2}) of $G$,
the right choice of function $\varkappa$ suitable to our purposes
arise from Remark \ref{rem4.1}. Indeed, due to formulae
(\ref{3.4}), if we define $\wtil G$ accordingly to (\ref{4.11})
then, for any $(x_0,x)\in\wtil G$, the angle $\theta_0$
between the $x_0$ axis and the
vector $|(x_0,x)|$ belongs to the interval
$(\phi,\pi-\phi)$. This leads at once to consider the following
cut-off function $\varkappa$:
\beqn\label{6.2}
&&\hskip -1,5truecm
\varkappa(x_0,x):={\cal E}\Big(\arccos\frac{x_0}{|(x_0,x)|}\Big)\,,\qq
(x_0,x)\in\wtil G\,,
\eeqn
where ${\cal E}\in C^{\infty}([0,\pi],\rsp_+)$,
$\|{\cal E}\|_{C([0,\pi])}=1$ and
for some $\ve\in (0,(\pi-2\phi)/6)$ satisfies
\beqn\label{6.3}
&&\hskip -1truecm
{\cal E}(\varphi)\equiv 0\,,\qq
\forall\,\varphi\in
[0,(\pi-6\ve)/2]\cup[(\pi+6\ve)/2,\pi]\,,
\\[1mm]
\label{6.4}
&&\hskip -1truecm
{\cal E}(\varphi)\equiv 1\,,
\,\qq\forall\,\varphi\in[(\pi-2\ve)/2,(\pi+2\ve)/2]\,.
\eeqn
In particular, by observing that our choice of $\ve$
guarantees ${\cal E}\equiv 1$ in an open interval containing $\pi/2$
and recalling $G=\{(x_0,x)\in\wtil G:x_0=0\}$,
we deduce that $\varkappa$ is equal to one on $G$
whereas ${\cal E}^{(l)}(\pi/2)=0$ for any $l\in\nsp\backslash\{0\}$.
Moreover, for any $l\in\nsp\cup\{0\}$ and $i=1,\ldots,n$ we have
\beqn\label{6.5}
&&\hskip -1,5truecm
D_{x_i}{\cal E}^{(l)}\Big(\arccos\frac{x_0}{|(x_0,x)|}\Big)=
\frac{x_0x_i}{|x||(x_0,x)|^2}\,{\cal E}^{(l+1)}
\Big(\arccos\frac{x_0}{|(x_0,x)|}\Big)\,,
\\[1mm]
\label{6.6}
&&\hskip -1,5truecm
D_{x_0}{\cal E}^{(l)}\Big(\arccos\frac{x_0}{|(x_0,x)|}\Big)=
-\frac{|x|}{|(x_0,x)|^2}\,{\cal E}^{(l+1)}
\Big(\arccos\frac{x_0}{|(x_0,x)|}\Big)\,.
\eeqn
Hence, when $\varkappa$ is defined by (\ref{6.2}),
from definition (\ref{6.1}) we derive
the following formulae for the first and the second derivatives of $v$,
where $i,j=1,\ldots,n$:
\beqn\label{6.7}
&&\hskip -0,8truecm
D_{x_i}v(x_0,x)=
\varkappa(x_0,x)e^{i\rho x_0}D_{x_i}u(x)
+[D_{x_i}\varkappa(x_0,x)]e^{i\rho x_0}u(x)\,,
\\[3mm]
\label{6.8}
&&\hskip -0,8truecm
D_{x_0}v(x_0,x)=
i\rho v(x_0,x)+[D_{x_0}\varkappa(x_0,x)]e^{i\rho x_0}u(x)\,,
\\[3mm]
\label{6.9}
&&\hskip -0,8truecm
D_{x_0}^2v(x_0,x)=
-\,\rho^2v(x_0,x)+\frac{|x|^2}{|(x_0,x)|^4}\,
{\cal E}''\Big(\arccos\frac{x_0}{|(x_0,x)|}\Big)e^{i\rho x_0}u(x)
\no
\\[2mm]
&&\hskip 1,8truecm
\,-2\Big[\frac{x_0-i\rho|(x_0,x)|^2}{|(x_0,x)|^2}\Big]
[D_{x_0}\varkappa(x_0,x)]e^{i\rho x_0}u(x)\,,
\\[3mm]
\label{6.10}
&&\hskip -0,8truecm
D_{x_i}D_{x_j}v(x_0,x)=
\varkappa(x_0,x)e^{i\rho x_0}D_{x_i}D_{x_j}u(x)
+2[D_{x_i}\varkappa(x_0,x)]e^{i\rho x_0}D_{x_j}u(x)
\no\\[2mm]
&&\hskip 2,4truecm
+\Big[\frac{\dl_{i,j}x_0}{|x||(x_0,x)|^2}-\frac{x_0x_ix_j(x_0^2+2|x|^2)}
{|x|^3|(x_0,x)|^4}\Big]
{\cal E}'\Big(\!\arccos\frac{x_0}{|(x_0,x)|}\Big)e^{i\rho x_0}u(x)
\no\\[2mm]
&&\hskip 2,4truecm
+\frac{x_0^2x_ix_j}{|x|^2|(x_0,x)|^4}
\,{\cal E}''\Big(\arccos\frac{x_0}{|(x_0,x)|}\Big)e^{i\rho x_0}u(x)\,,
\\[3mm]
\label{6.11}
&&\hskip -0,8truecm
D_{x_0}D_{x_j}v(x_0,x)=
[i\rho\,\varkappa(x_0,x)+D_{x_0}\varkappa(x_0,x)]
e^{i\rho x_0}D_{x_j}u(x)
\no\\[2mm]
&&\hskip 2,4truecm
+\,\Big[\frac{i\rho x_0x_j}{|x||(x_0,x)|^2}
+\frac{x_j(|x|^2-x_0^2)}{|x||(x_0,x)|^4}\Big]
{\cal E}'\Big(\!\arccos\frac{x_0}{|(x_0,x)|}\Big)e^{i\rho x_0}u(x)
\no\\[2mm]
&&\hskip 2,4truecm
-\,\frac{x_0x_j|x|}{|x||(x_0,x)|^4}\,
{\cal E}''\Big(\!\arccos\frac{x_0}{|(x_0,x)|}\Big)
e^{i\rho x_0}u(x)\,.
\eeqn
In addition,  using (\ref{6.5}), (\ref{6.6}) with $l=0$ we deduce,
for any $j,k=1,\ldots,n$,
\beqn\label{6.12}
&&\hskip -1,5truecm
D_{x_j}D_{x_k}\varkappa(x_0,x)
=\frac{x_0^2x_jx_k}{|x|^2|(x_0,x)|^4}\,
{\cal E}''\Big(\arccos\frac{x_0}{|(x_0,x)|}\Big)\no
\\[2mm]
&&\hskip 1,9truecm
-\,\frac{\dl_{j,k}x_0D_{x_0}\varkappa(x_0,x)}{|x|^2}
-\frac{x_k(x_0^2+3|x|^2)D_{x_j}\varkappa(x_0,x)}{|x|^2|(x_0,x)|^2}\,.
\eeqn
Now, let ${\cal A}_\psi(x;D_x,D_{x_0})$
and ${\cal B}((x_0,x);D_x,D_{x_0})$ be defined as in
the statement of Lemma \ref{lem4.2}.
Through easy but lengthy computations,
from (\ref{6.7}), (\ref{6.9}), (\ref{6.12}) we obtain
\beqn\label{6.13}
&&\hskip -0,5truecm
{\cal A}_\psi(x;D_x,D_{x_0})v(x_0,x)\no
\\
&&\hskip -1truecm
=\varkappa(x_0,x)e^{i\rho x_0}
[{\cal A}(x;D_x)-\rho^2e^{i\psi}I]u(x)
+2e^{i\rho x_0}
\sum_{j,k=1}^na_{j,k}(x)D_{x_k}u(x)D_{x_j}\varkappa(x_0,x)\no
\\
&&\hskip -0,5truecm
+\,{\cal E}''\Big(\arccos\frac{x_0}{|(x_0,x)|}\Big)e^{i\rho x_0}u(x)
\sum_{j,k=1}^n\frac{a_{j,k}(x)x_0^2x_jx_k}{|x|^2|(x_0,x)|^4}\no
\\
&&\hskip -0,5truecm
+e^{i(\psi+\rho x_0)}u(x)
\bigg\{\frac{|x|^2}{|(x_0,x)|^4}\,
{\cal E}''\Big(\arccos\frac{x_0}{|(x_0,x)|}\Big)
-2\Big[\frac{x_0-i\rho|(x_0,x)|^2}{|(x_0,x)|^2}\Big]
D_{x_0}\varkappa(x_0,x)\bigg\}\no
\\[1mm]
&&\hskip -0,5truecm
+\,e^{i\rho x_0}u(x)
\bigg\{
\sum_{j,k=1}^nD_{x_k}a_{j,k}(x)D_{x_j}\varkappa(x_0,x)
-\sum_{j=1}^n\frac{a_{j,j}(x)x_0D_{x_0}\varkappa(x_0,x)}{|x|^2}\no
\\
&&\hskip 1,8truecm
-\sum_{j,k=1}^n\frac{a_{j,k}(x)x_k(x_0^2+3|x|^2)
D_{x_j}\varkappa(x_0,x)}{|x|^2|(x_0,x)|^2}
\bigg\}\no
\\
&&\hskip -1truecm
=:\sum_{l=1}^5J_l(u,{\cal E},(x_0,x))\,,
\eeqn
whereas, from (\ref{6.7}) and (\ref{6.8}), we get
\beqn\label{6.14}
&&\hskip -1,5truecm
{\cal B}((x_0,x);D_x,D_{x_0})v(x_0,x)=
\varkappa(x_0,x){e^{i\rho x_0}}{\cal B}(x;D_x)u(x)
+i\rho x_0\varkappa(x_0,x)e^{i\rho x_0}u(x)
\no\\[1mm]
&&\hskip 3,8truecm
+\,e^{i\rho x_0}u(x)\Big[x_0D_{x_0}\varkappa(x_0,x)+\sum_{i=1}^n
b_i(x)D_{x_i}\varkappa(x_0,x)\Big]
\no\\
&&\hskip 3,3truecm
=:\sum_{l=6}^8J_l(u,{\cal E},(x_0,x))\,.
\eeqn
In the next section,
with the help of (\ref{6.13}) and (\ref{6.14}),
we will upper bound
the norms of ${\cal A}_\psi(x;D_x,D_{x_0})v$
and ${\cal B}((x_0,x);D_x,D_{x_0})v$, respectively in $L_{p,-1}(\wtil G)$
and $W_{p,-1}^{1-p^{-1}}(\partial \wtil G)$,
in terms of the $W_{p,-1}^k(G)$-norms , $k=0,1,2$, of $u$
and of an any extension to $G$ of its assigned boundary values.
Just these estimates will be the argument of the forthcoming
lemmata Lemma \ref{lem7.3} and Lemma \ref{lem7.4},
which will be a fundamental step in the proof of our main result
Theorem \ref{thm4.6}.
\section{Proof of Theorem \ref{thm4.6}}
\setcounter{equation}{0}
As we said at the end of Section 6, Theorem \ref{thm4.6} will be an
easy consequence of two crucial lemmata, Lemma \ref{lem7.3}
and Lemma \ref{lem7.4}. We postpone such lemmata
to the following considerations which strictly depends
on the class (\ref{4.2}) of domains $G$ we restrict to work with.
First, observe that $|(x_0,x)|\ge |x|$ implies
\beqn\label{7.4}
&&\hskip -2,5truecm
|(x_0,x)|^{-q}
\le
|x|^{-q}\,,\qq
\forall\,q\ge 1\,.
\eeqn
Hence, if $w(x_0,x)$ is a function such that
$|w(x_0,x)|\le|(x_0,x)|^{-k}|w_1(x_0,x)||w_2(x)|$,
for some $k\in\nsp\cup\{0\}$ and some $w_1\in C(\wtil G)$
vanishing for $x_0\notin(-\dl_0,\dl_0)$, $0<\dl_0<\|\eta\|_{C([0,C_2])}$,
then from (\ref{7.4})
we easily find
\beqn\label{7.5}
&&\hskip -1truecm
\|w\|_{L_{p,-1}(\wtil G)}\le
(2\dl_0)^{1/p}\|w_1\|_{C(\wtil G)}
\|w_2\|_{L_{p,-(k+1)}(G)}\,.
\eeqn
Moreover, assumptions i), ii) on function
$\eta$ which describes the boundary $\partial G$ of $G$ imply
\beqn\label{7.6}
&&\hskip -1truecm
|(x_0,x)|\le\{[\eta(x_n)]^2+x_n^2\}^{1/2}
\le\{C_3^2+1\}^{1/2}x_n\,,
\q\;\forall\,(x_0,x)\in\wtil G\,.
\eeqn
Therefore, if we set $C_4=\{C_3^2+1\}^{1/2}$
and we use $x_n\le |x|$, from (\ref{7.6}) we deduce
\beqn\label{7.7}
C_4^{-q}|x|^{-q}\le|(x_0,x)|^{-q}\,,\qq
\forall\,q\ge 1\,,\qq\forall\,(x_0,x)\in\wtil G\,,
\eeqn
and hence, for any $w\in L_{p,-1}(\wtil G)$,
we deduce also
$\|w\|_{L_{p,-1}(\wtil G)}\ge
C_4^{-1}\|w\|_{L_{p,-1}(G)}$.
\begin{lemma}\label{lem7.3}
Let $p>n$ and $u\in W_{p,-1}^2(G)$, where $G$ is defined by (\ref{4.2}),
and let ${\cal A}(x;D_x)$ be the differential
operator (\ref{1.1}) with coefficients
$a_{i,j}$, $i,j=1,\ldots,n$, satisfying (\ref{4.3}).
Then, when $v$ is defined by (\ref{6.1}), (\ref{6.2})
and ${\cal A}_\psi(x;D_x,D_{x_0})$ is defined as in the statement
of Lemma \ref{lem4.2}, for any $\rho>0$ and
some $\dl_0(\ve)>0$  the following estimate holds:
\beqn\label{7.9}
&&
\|{\cal A}_\psi(x;D_x,D_{x_0})v\|_{L_{p,-1}(\wtil G)}
\no
\\[2mm]
&&\hskip -1truecm
\le [2\dl_0(\ve)]^{1/p}
\big\{\|\big({\cal A}(x;D_x)-\rho^2e^{i\psi} I\big)u\|_{L_{p,-1}(G)}
+M_1(1+\rho)\|u\|_{W_{p,-1}^2(G)}\big\}\,,
\eeqn
The positive constant $M_1$ in (\ref{7.9}) depends only on $p$, $n$,
the $C^1(G)$-norm of the coefficients of ${\cal A}(x;D_x)$
and the constants $C_j$, $j=2,3$,
intervening in the properties i)--iii)
for the function $\eta$ which describes the boundary $\partial G$
of $G$.
\end{lemma}
\begin{proof}
Since from formula (\ref{6.13}) it follows
\beqn\label{7.10}
\|{\cal A}_\psi(x;D_x,D_{x_0})v\|_{L_{p,-1}(\wtil G)}
\le\sum_{l=1}^5\|J_l(u,{\cal E},(x_0,x))\|_{L_{p,-1}(\wtil G)}\,.
\eeqn
we need only to estimate from above each norm
$\|J_l(u,{\cal E},(x_0,x))\|_{L_{p,-1}(\wtil G)}$, $l=1,\ldots,5$,
and then to rearrange the term.
First, from (\ref{7.5}) and $\|{\cal E}\|_{C([0,\pi])}= 1$
we immediately get
\beqn\label{7.11}
&&\hskip -2 truecm
\|J_1(u,{\cal E},(x_0,x))\|_{L_{p,-1}(\wtil G)}
\le  [2\dl_0(\ve)]^{1/p}
\|\big({\cal A}(x;D_x)-\rho^2e^{i\psi}I\big)u\|_{L_{p,-1}(G)}\,.
\eeqn
where using (\ref{6.3}), (\ref{6.4}) and (\ref{7.6}) we have
set
\beqn\label{7.12}
&&\hskip -1truecm
\dl_0(\ve)=C_4C_2\cos[(\pi-6\ve)/2],\qq\ve\in (0,(\pi-2\phi)/6)\,.
\eeqn
Now, observe that
for any $(x_0,x)\in\wtil G$ and $l\in\nsp\cup\{0\}$
from (\ref{6.5}), (\ref{6.6}) we derive
\beqn\label{7.13}
&&\hskip -1truecm
\big|D_{x_j}{\cal E}^{(l)}\Big(\arccos\frac{x_0}{|(x_0,x)|}\Big)\big|\le
\frac{1}{|(x_0,x)|}\,\|{\cal E}^{(l+1)}\|_{C([0,\pi])}\,,\q
\,\q j=0,\ldots,n\,,
\eeqn
whereas, since $u\in W_{p,-1}^2(G)$ and $-1<-np^{-1}$,
from Theorem \ref{thm2.1} it follows
\beqn\label{7.14}
&&\hskip -1,5truecm
\big\{\|u\|_{L_{p,-3}(G)}^p
+\sum_{k=1}^n\|D_{x_k}u\|_{L_{p,-2}(G)}^p\big\}^{1/p}\le
\|u\|_{V_{p,-1}^2(G)}\le
c\|u\|_{ W_{p,-1}^2(G)}\,.
\eeqn
Hence, recalling the definition of $J_2(u,{\cal E},(x_0,x))$ in
(\ref{6.13}) and that of $\varkappa$ in (\ref{6.2}),
from (\ref{7.5}), (\ref{7.13}) and (\ref{7.14}) we deduce
\beqn\label{7.15}
&&\hskip -0,5truecm
\|J_2(u,{\cal E},(x_0,x))\|_{L_{p,-1}(\wtil G)}\no
\\
&&\hskip -1truecm
\le 2n
[2\dl_0(\ve)]^{1/p}
\|{\cal E}'\|_{C([0,\pi])}
\max_{j,k=1,\ldots,n}\|a_{j,k}\|_{C(G)}
\Big(\sum_{k=1}^n\|D_{x_k}u\|_{L_{p,-2}(G)}^p\Big)^{1/p}
\no
\\
&&\hskip -1truecm
\le 2cn [2\dl_0(\ve)]^{1/p}\|{\cal E}'\|_{C([0,\pi])}
\max_{j,k=1,\ldots,n}\|a_{j,k}\|_{C(G)}\|u\|_{W_{p,-1}^2(G)}\,,
\eeqn
and similarly, but taking advantage from
\beqn\label{7.16}
&&\hskip -1,5truecm
\Big|\frac{x_0^2x_jx_k}{|x|^2|(x_0,x)|^4}\Big|\le
\frac{1}{|(x_0,x)|^2}\,,\qq
j,k=1,\ldots,n\,,
\eeqn
we obtain
\beqn\label{7.17}
&&\hskip -1,8truecm
\|J_3(u,{\cal E},(x_0,x))\|_{L_{p,-1}(\wtil G)}\le
n^2[2\dl_0(\ve)]^{1/p}\|{\cal E}''\|_{C([0,\pi])}
\max_{j,k=1,\ldots,n}\|a_{j,k}\|_{C(G)}\|u\|_{L_{p,-3}(G)}
\no
\\[1mm]
&&\hskip 2,7truecm
\le
cn^2[2\dl_0(\ve)]^{1/p}
\|{\cal E}''\|_{C([0,\pi])}\max_{j,k=1,\ldots,n}\|a_{j,k}\|_{C(G)}
\|u\|_{W_{p,-1}^2(G)}.
\eeqn
Finally, using (\ref{7.6}) and (\ref{7.13}), it is easy
to prove that the factors
on the braces in the definition of
$J_4(u,{\cal E},(x_0,x))$ and $J_5(u,{\cal E},(x_0,x))$
have their absolute values which are bounded from above respectively
by $|(x_0,x)|^{-2}(3+2\rho C_4C_2)\|{\cal E}\|_{C^2([0,\pi])}$
and  $[|x||(x_0,x)|]^{-1}n^2[4+C_4(1+C_2)]
\|{\cal E}'\|_{C([0,\pi])}\max_{j,k=1,\ldots,n}\|a_{j,k}\|_{C^1(G)}$.\\
Therefore, if we set $M_0=\{\max[2^{p-1}3^p,2^{2p-1}C_4^pC_2^p]\}^{1/p}$,
from (\ref{7.5}) and (\ref{7.14}) we find
\beqn\label{7.18}
&&\hskip -1,5truecm
\|J_4(u,{\cal E},(x_0,x))\|_{L_{p,-1}(\wtil G)}
\le [2\dl_0(\ve)]^{1/p}M_0(1+\rho)
\|{\cal E}\|_{C^2([0,\pi])}\|u\|_{L_{p,-3}(G)}
\no\\
&&\hskip 3truecm
\le c[2\dl_0(\ve)]^{1/p}M_0
(1+\rho)\|{\cal E}\|_{C^2([0,\pi])}
\|u\|_{W_{p,-1}^2(G)}\,,
\eeqn
\beqn\label{7.19}
&&\hskip -1,5truecm
\|J_5(u,{\cal E},(x_0,x))\|_{L_{p,-1}(\wtil G)}\no
\\[1mm]
&&\hskip -2truecm
\le n^2[2\dl_0(\ve)]^{1/p}[4+C_4(1+C_2)]
\|{\cal E}'\|_{C([0,\pi])}\max_{j,k=1,\ldots,n}\|a_{j,k}\|_{C^1(G)}
\|u\|_{L_{p,-3}(G)}\,
\no
\\[1mm]
&&\hskip -2truecm
\le
cn^2[2\dl_0(\ve)]^{1/p}[4+C_4(1+C_2)]\|{\cal E}'\|_{C([0,\pi])}
\max_{j,k=1,\ldots,n}\|a_{j,k}\|_{C^1(G)}
\|u\|_{W_{p,-1}^2(G)}\,.
\eeqn
By replacing (\ref{7.11}), (\ref{7.15}),
(\ref{7.17})--(\ref{7.19}) in (\ref{7.10})
and rearranging the term we obtain (\ref{7.9}) with
$M_1=c\{M_0+n[2+5n+nC_4(1+C_2)]
\max_{j,k=1,\ldots,n}\|a_{j,k}\|_{C^1(G)}\}
\|{\cal E}\|_{C^2([0,\pi])}$.
\end{proof}
\begin{lemma}\label{lem7.4}
Let $p>n$ and $u\in W_{p,-1}^2(G)$ where $G$ is defined by (\ref{4.2})
and let ${\cal B}(x;D_x)$ be the differential
operator (\ref{4.6}) with coefficients
$b_j$, $j=0,\ldots,n$, satisfying (\ref{4.7}).
Then, $g_0\in W_{p,-1}^1(G)$ being any extension
to $G$ of ${\cal B}(x;D_x)u$, when $v$ is defined by (\ref{6.1}), (\ref{6.2})
and ${\cal B}((x_0,x);D_x,D_{x_0})$ is defined as in the statement
of Lemma \ref{lem4.2}, for any $\rho>0$
and some $\dl_0(\ve)>0$ the following estimate holds:
\beqn\label{7.20}
&&\hskip -1truecm
\|{\cal B}((x_0,x);D_x,D_{x_0})v\|_{W_{p,-1}^{1-p^{-1}}(\partial\wtil G)}
\no\\
&&\hskip -1,5truecm
\le
2[\dl_0(\ve)]^{1/p}M_2\big\{
\|g_0\|_{W_{p,-1}^1(G)}+\rho\|g_0\|_{L_{p,-1}(G)}
+\|u\|_{W_{p,-1}^2(G)} \no
\\[1mm]
&&\hskip 1,7truecm
+\,(1+2\rho)\|u\|_{W_{p,-1}^1(G)}+(\rho+\rho^2)\|u\|_{L_{p,-1}(G)}\big\}\,,
\eeqn
The constant $M_2>1$ in (\ref{7.20}) depends only on $p$, $n$,
the $C^1(G)$-norm of the coefficients of ${\cal B}(x;D_x)$
and the constants $C_j$, $j=2,3$,
intervening in the properties i)--iii)
for the function $\eta$ which describes the boundary $\partial G$
of $G$.
\end{lemma}
\begin{proof}
First, from (\ref{6.14}) we get
\beqn\label{7.21}
&&\hskip -1,1truecm
\|{\cal B}((x_0,x);D_x,D_{x_0})v\|_{W_{p,-1}^{1-p^{-1}}(\partial\wtil G)}\le
\sum_{l=6}^8
\|J_l(u,{\cal E},(x_0,x))\|_{W_{p,-1}^{1-p^{-1}}(\partial\wtil G)}\,,
\eeqn
and observe that, due to the definition (\ref{4.11}) of $\wtil G$,
if $(x_0,x)$ belong to $\partial\wtil G$ then it is
of the form $(0,x)$ with $x\in \partial G$
or $(x_0,x)$ with $x_0\neq 0$ and $x\in G$.
Therefore (cf. also (\ref{4.7})), the term $J_l(u,{\cal E},(x_0,x))$,
$l=6,7,8$, in (\ref{7.21}) are well defined for any
$(x_0,x)\in \wtil G\cup\partial\wtil G$.
Hence, recalling the definition (\ref{2.3})
of the norm in the spaces of traces and
using (\ref{6.7}), (\ref{6.8}) with $u$ replaced
by $g_0$ and the inequality $|a+b|^q\le 2^{q-1}(|a|^q+|b|^q)$,
$a,b\in\csp$, $q\ge 1$, from (\ref{6.5}), (\ref{6.6}),
(\ref{7.13}) and (\ref{7.5})
 we obtain
\beqn\label{7.22}
&&\hskip -1truecm
\|J_6(u,{\cal E},(x_0,x))\|_{W_{p,-1}^{1-p^{-1}}(\partial\wtil G)}
\le \|J_6(u,{\cal E},(x_0,x))\|_{W_{p,-1}^1(\wtil G)}\,\no
\\
&&\hskip -1,5truecm
\le 2[\dl_0(\ve)]^{1/p}\big\{
\|g_0\|_{W_{p,-1}^1(G)}
+(n+1)^{1/p}\|{\cal E}'\|_{C([0,\pi])}\|g_0\|_{L_{p,-2}(G)}
+\rho\|g_0\|_{L_{p,-1}(G)}\big\},
\eeqn
$\dl_0(\ve)$ being defined by (\ref{7.12}).
Now, for any $w\in W_{p,-1}^1(G)$ with $p>n$ Theorem \ref{thm2.1} imply
\beqn\label{7.23}
&&\hskip -1truecm
\big\{\|w\|_{L_{p,-2}(G)}^p
+\sum_{k=1}^n\|D_{x_k} w\|_{L_{p,-1}(G)}^p\big\}^{1/p}
\le
\|w\|_{V_{p,-1}^1(G)}\le
c\|w\|_{ W_{p,-1}^1(G)}\,,
\eeqn
and consequently, if we set $M_3=[1+c(n+1)^{1/p}\|{\cal E}'\|_{C([0,\pi])}]$,
from (\ref{7.22}) we get
\beqn\label{7.24}
&&\hskip -2truecm
\|J_6(u,{\cal E},(x_0,x))\|_{W_{p,-1}^{1-p^{-1}}(\partial\wtil G)}\le
2[\dl_0(\ve)]^{1/p}\big\{
M_3\|g_0\|_{W_{p,-1}^1(G)}+\rho\|g_0\|_{L_{p,-1}(G)}\big\}\,.
\eeqn
Similarly, from (\ref{6.5})--(\ref{6.8}),
(\ref{7.13}), (\ref{7.5}) and (\ref{7.23}) we obtain
\beqn\label{7.25}
&&\hskip -1truecm
\|J_7(u,{\cal E},(x_0,x))\|_{W_{p,-1}^{1-p^{-1}}(\partial\wtil G)}
\le \|J_7(u,{\cal E},(x_0,x))\|_{W_{p,-1}^1(\wtil G)}\,\no
\\[1mm]
&&\hskip -1,5truecm
\le 2M_4[\dl_0(\ve)]^{1/p}
\big\{ M_3\rho\|u\|_{W_{p,-1}^1(G)}
+2^{1-1/p}(\rho+\rho^2)\|u\|_{L_{p,-1}(G)}\big\}\,,
\eeqn
where $M_4=\max[1,\dl_0(\ve)]$.\\
Before to estimate the term $J_8(u,{\cal E},(x_0,x))$
in (\ref{7.21}) observe that from (\ref{6.5}) and (\ref{6.6})
it follows
\beqn\label{7.26}
&&\hskip -1truecm
x_0D_{x_0}\varkappa(x_0,x)+\sum_{i=1}^n
b_i(x)D_{x_i}\varkappa(x_0,x)
\no\\
&&\hskip -1,5truecm
=-\Big[\frac{x_0|x|}{|(x_0,x)|^2}-\sum_{i=1}^n
\frac{b_i(x)x_0x_i}{|x||(x_0,x)|^2}\Big]
{\cal E}'\Big(\arccos\frac{x_0}{|(x_0,x)|}\Big)\,,
\eeqn
so that, using (\ref{7.6}), we easily get
\beqn\label{7.27}
&&\hskip -1,5truecm
\Big|x_0D_{x_0}\varkappa(x_0,x)+\sum_{i=1}^n
b_i(x)D_{x_i}\varkappa(x_0,x)
\Big|\le
\frac{M_5\|{\cal E}'\|_{C([0,\pi])}}{|(x_0,x)|}\,,
\eeqn
where $M_5=[n\max_{i=1,\ldots,n}\|b_i\|_{C(G)}+C_4C_2]$.
In addition,
for any $k=1,\ldots,n$ we have
\beqn
&&\hskip -1,5truecm
D_{x_0}\Big[\frac{x_0|x|}{|(x_0,x)|^2}-
\sum_{i=1}^n\frac{b_i(x)x_0x_i}{|x||(x_0,x)|^2}\Big]
=\frac{(|x|^2-x_0^2)}{|(x_0,x)|^4}\Big[|x|-\sum_{i=1}^n
\frac{b_i(x)x_i}{|x|}\Big],\no
\\[1mm]
&&\hskip -1,5truecm
D_{x_k}\Big[\frac{x_0|x|}{|(x_0,x)|^2}-
\sum_{i=1}^n\frac{b_i(x)x_0x_i}{|x||(x_0,x)|^2}\Big]
\no\\[1mm]
&&\hskip -1,5truecm
=\frac{x_0x_k(x_0^2-|x|^2)}{|x||(x_0,x)|^4}
+\sum_{i=1}^n\frac{x_0[x_iD_{x_k}b_i(x)+\dl_{i,k}b_i(x)]}{|x||(x_0,x)|^2}
-\sum_{i=1}^n\frac{b_i(x)x_0x_ix_k[x_0^2+3|x|^2]}{|x|^3|(x_0,x)|^4}\,,
\no
\eeqn
and hence, applying the Leibniz's formula to the right-hand
side of (\ref{7.26}) and using (\ref{7.4}), (\ref{7.6}), (\ref{7.13}) and
(\ref{7.16}), it is easy to obtain
\beqn\label{7.28}
&&\hskip -1truecm
\Big|D_{x_0}\Big[x_0D_{x_0}\varkappa(x_0,x)+\sum_{i=1}^n
b_i(x)D_{x_i}\varkappa(x_0,x)\Big]\Big|
\le \frac{M_5}
{|(x_0,x)|^2}\sum_{k=1}^2\|{\cal E}^{(k)}\|_{C([0,\pi])}\,,
\\[1mm]
\label{7.29}
&&\hskip -1truecm
\Big|D_{x_k}\Big[x_0D_{x_0}\varkappa(x_0,x)+\sum_{i=1}^n
b_i(x)D_{x_i}\varkappa(x_0,x)\Big]\Big|
\le\frac{M_6}{|x||(x_0,x)|}\sum_{k=1}^2\|{\cal E}^{(k)}\|_{C([0,\pi])}\,,
\eeqn
where $M_6=[n(5+C_4C_2)\max_{i=1,\ldots,n}\|b_i\|_{C^1(G)}+C_4C_2]$
and in (\ref{7.29}) we have used $M_5<M_6$.\\
Therefore,
combining (\ref{7.5}) and (\ref{7.27})--(\ref{7.29}) we deduce
\beqn\no
&&\hskip -1truecm
\|J_8(u,{\cal E},(x_0,x))\|_{W_{p,-1}^{1-p^{-1}}(\partial\wtil G)}
\le \|J_8(u,{\cal E},(x_0,x))\|_{W_{p,-1}^1(\wtil G)}\,\no
\\[1mm]
&&\hskip -1,5truecm
\le 2M_6[\dl_0(\ve)]^{1/p}\|{\cal E}\|_{C^2([0,\pi])}
\big\{2^{1-1/p}(1+\rho)\|u\|_{L_{p,-2}(G)}+2^{1-1/p}
(1+2n)^{1/p}\|u\|_{L_{p,-3}(G)}\no
\\
&&\hskip 3,6truecm
+\big[\|u\|_{L_{p,-3}(G)}^p
+\sum_{k=1}^n\|D_{x_k}u\|_{L_{p,-2}(G)}^p\big]^{1/p}\,\big\}\,,
\no
\eeqn
and hence, using (\ref{7.14}) and (\ref{7.23}),
\beqn\label{7.30}
&&\hskip -1truecm
\|J_8(u,{\cal E},(x_0,x))\|_{W_{p,-1}^{1-p^{-1}}(\partial\wtil G)}
\no\\
&&\hskip -1,5truecm
\le 2M_7[\dl_0(\ve)]^{1/p}
\big\{M_8\|u\|_{W_{p,-1}^2(G)} +2^{1-1/p}(1+\rho)\|u\|_{W_{p,-1}^1(G)}\big\}\,,
\eeqn
where $M_7=cM_6\|{\cal E}\|_{C^2([0,\pi])}$
and $M_8=[1+2^{1-1/p}(1+2n)^{1/p}]$.
Rearranging (\ref{7.24}), (\ref{7.25}) and (\ref{7.30})
from (\ref{7.21}) we derive (\ref{7.20}) with the constant
$M_2=\max[M_4,M_7]\times\max[M_3,M_8]$.
\end{proof}
We can now prove the main result of the paper.
To simplify notations, from now on for any $u\in W_{p,-1}^l(G)$, $l\ge 0$,
we will set $\|D^lu\|_{L_{p,-1}(G)}=\sum_{|\a|=l}\|D^{\a}u\|_{L_{p,-1}(G)}$.\\
{\it Proof of Theorem \ref{thm4.6}.}
For every $u\in W_{p,-1}^2(G)$
and $\rho>0$ we define the function $v$ accordingly
to (\ref{6.1}) (\ref{6.2})
and we observe that $p>n$, $p\neq n+1$, imply
\beqn\no
&&\hskip -1truecm
\left\{\!
\begin{array}{lll}
-1<-(n+1)p^{-1}\,,\qq{\rm if}\q p>n+1\,,
\\[2mm]
-(n+1)p^{-1}<-1<1-(n+1)p^{-1}\,, \qq{\rm if}\q n<p<n+1\,.
\end{array}
\right.
\eeqn
Therefore, recalling (\ref{3.22}),
the assumptions of Theorem \ref{thm3.2} for $\b=-1$ are both
satisfied, the second one with $\nu=0$.
Moreover, since we have assumed problem (\ref{4.14}) to be regular
in the sense of Definition \ref{def4.5},
we can apply to the function $v$ the estimate (\ref{3.21})
with the choice of the parameter $l$, $\vec t$, $\vec s$ and $\vec\s$
as in (\ref{3.22}):
\beqn\label{7.32}
&&\hskip -2truecm
\|v\|_{W_{p,-1}^{2}(\wtil G)}\le
c_1\big\{
\|{\cal A}_\psi(x;D_x,D_{x_0})v\|_{L_{p,-1}(\wtil G)}
\no\\[1mm]
&&\hskip 1,1truecm
+\|{\cal B}((x_0,x);D_x,D_{x_0})v\|_{W_{p,-1}^{1-p^{-1}}(\partial\wtil G)}
+\|v\|_{W_{p,-1}^{1}(\wtil G)}\big\}.
\eeqn
Since the norms $L_{p,-1}(\wtil G)$ and 
$W_{p,-1}^{1-p^{-1}}(\partial\wtil G)$ of
${\cal A}_\psi(x;D_x,D_{x_0})v$ and ${\cal B}((x_0,x);D_x,D_{x_0})v$
have been estimated in Lemma \ref{lem7.3} and Lemma \ref{lem7.4}, 
respectively, it remains only to analyze the term 
$\|v\|_{W_{p,-1}^1(\wtil G)}$ in (\ref{7.32}). 
But, due to definition (\ref{6.1}),
as in (\ref{7.22}) and (\ref{7.24}) with $g_0$ replaced by $u$,
we obtain
\beqn\label{7.33}
&&\hskip -0,7truecm
\|v\|_{W_{p,-1}^1(\wtil G)}
\le 2[\dl_0(\ve)]^{1/p}\big\{
M_3\|u\|_{W_{p,-1}^1(G)}+\rho\|u\|_{L_{p,-1}(G)}\big\}\,.
\eeqn
Since for any $w\in W_{p,-1}^l(G)$, $l\ge 0$,
we have $\|w\|_{W_{p,-1}^l(G)}\le \sum_{|\a|=0}^l\|D^{\a}w\|_{L_{p,-1}(G)}$,
if we set $M_9=M_2\max[1,M_1]$,
$M_2$ being defined at the end of the proof of Lemma \ref{lem7.4},
by combining (\ref{7.9}), (\ref{7.20}) and (\ref{7.33})
from (\ref{7.32}) we obtain
\beqn\label{7.34}
&&\hskip -0,8truecm
\|v\|_{W_{p,-1}^2(\wtil G)}
\le 2c_1M_9[\dl_0(\ve)]^{1/p}\Big\{
\|({\cal A}(x;D_x)-\rho^2e^{i\psi}I)u\|_{L_{p,-1}(G)}
+(4+5\rho+\rho^2)\|u\|_{L_{p,-1}(G)}
\no\\[1mm]
&&\hskip 4,8truecm
+(4+3\rho)\|Du\|_{L_{p,-1}(G)}
+(2+\rho)\|D^2u\|_{L_{p,-1}(G)}
\no\\[1mm]
&&\hskip 4,8truecm
+\,(1+\rho)\|g_0\|_{L_{p,-1}(G)}+\|Dg_0\|_{L_{p,-1}(G)}\Big\}\,,
\eeqn
On the other hand, using $G=\{(x_0,x)\in\wtil G:x_0=0\}$, (\ref{6.4})  and
${\cal E}^{(k)}(\pi/2)=0$, $k=1,2$,
from (\ref{6.1}) and (\ref{6.7})--(\ref{6.11})
for any $(x_0,x)\in\wtil G$ we deduce the following inequalities
\beqn\no
&&\hskip -1truecm
\left\{
\begin{array}{lll}
|v(x_0,x)|\ge|v(0,x)|=|u(x)|\,,
\\[2mm]
|D_{x_i}v(x_0,x)|\ge|D_{x_i}v(0,x)|=|D_{x_i}u(x)|\,,\q
i=1,\ldots,n\,,
\\[2mm]
|D_{x_0}v(x_0,x)|\ge|D_{x_0}v(0,x)|=\rho|u(x)|\,,
\\[2mm]
|D_{x_0}^2v(x_0,x)|\ge|D_{x_0}^2v(0,x)| =\rho^2|u(x)|\,,
\\[2mm]
|D_{x_i}D_{x_j}v(x_0,x)|
\ge|D_{x_i}D_{x_j}v(0,x)| =|D_{x_i}D_{x_j}u(x)|\,,
\q i,j=1,\ldots,n,
\\[2mm]
|D_{x_0}D_{x_j}v(x_0,x)|
\ge|D_{x_0}D_{x_j}v(0,x)|=\rho|D_{x_j}u(x)|\,,\q j=1,\ldots,n.
\end{array}
\right.
\eeqn
Hence, using (\ref{7.7}) we obtain
\beqn
&&\hskip -0,4truecm
C_4^p\,\|v\|_{W_{p,-1}^2(\wtil G)}^p
=\sum_{0\le|\a|\le2}
\int_{\wtil G}C_4^p|(x_0,x)|^{-p}\big|D^{\a}v(x_0,x)\big|^p\de x_0\de x
\no\\
&&\hskip -0,8truecm
\ge\int_{\wtil G}|x|^{-p}\Big[
|v(x_0,x)|^p+\sum_{i=0}^n|D_{x_i}v(x_0,x)|^p
+\sum_{i,j=0}^n|D_{x_i}D_{x_i}v(x_0,x)|^p\Big]\de x_0\de x
\no\\
&&\hskip -0,8truecm
\ge\int_{G}|x|^{-p}\Big[
(1+\rho^p+\rho^{2p})|u(x)|^p+(1+\rho^p)\sum_{i=1}^n|D_{x_i}u(x)|^p
+\sum_{i,j=1}^n|D_{x_i}D_{x_i}u(x)|^p\Big]\de x
\no\\[1mm]
&&\hskip -0,8truecm
\ge \rho^{2p}\|u\|_{L_{p,-1}(G)}^p
+\rho^p\|Du\|_{L_{p,-1}(G)}^p
+\|D^2u\|_{L_{p,-1}(G)}^p\,.
\no
\eeqn
Taking into account (\ref{7.34}), it follows
\beqn\label{7.35}
&&\hskip -0,5truecm
\rho^{2}\|u\|_{L_{p,-1}(G)}
+\rho
\|Du\|_{L_{p,-1}(G)}
+\|D^2u\|_{L_{p,-1}(G)}\le
3C_4\,\|v\|_{W_{p,-1}^2(\wtil G)}
\no\\[1mm]
&&\hskip -1truecm
\le
M_{10}(\ve)\Big\{
\|\big({\cal A}(x;D_x)-\rho^2e^{i\psi}I\big)u\|_{L_{p,-1}(G)}
+(4+5\rho+\rho^2)\|u\|_{L_{p,-1}(G)}
\no\\[1mm]
&&\hskip 1truecm
+\,(4+3\rho)\|Du\|_{L_{p,-1}(G)}+(2+\rho)\|D^2u\|_{L_{p,-1}(G)}
\no\\[1mm]
&&\hskip 1truecm
+\,(1+\rho)\|g_0\|_{L_{p,-1}(G)}+\|Dg_0\|_{L_{p,-1}(G)}\Big\}
\eeqn
where we have set $M_{10}(\ve)=6c_1C_4M_9[\dl_0(\ve)]^{1/p}$.
Now, from (\ref{7.12}) we deduce that
$M_{10}(\ve)$ goes to zero as $\ve\to 0^+$.
Therefore, if we take $\l=\rho^2e^{i\psi}$
and we assume $\ve$ sufficiently small,
we can take $\rho$ so large so that the following inequalities
are satisfied
\beqn\label{7.36}
\left\{
\begin{array}{lll}
M_{10}(\ve)(4+5\rho+\rho^2)\le \rho^2/2\,,
\\[2mm]
M_{10}(\ve)(4+3\rho)\le \rho/2\,,
\\[2mm]
M_{10}(\ve)(2+\rho)\le 1/2\,.
\end{array}
\right.
\eeqn
From (\ref{7.35}) and (\ref{7.36}) our statement follows
with $M=2M_{10}(\ve)$ in (\ref{4.16}).
\begin{remark}\label{rem7.6}
\emph{
In the latter part of the proof of Theorem \ref{thm4.6} 
we have assumed $\ve$ to be close to zero which, equivalently,
means that the function ${\cal E}$ in (\ref{6.2})
has its support in a small neighborhood of $\pi/2$
(cf. (\ref{6.3}) and (\ref{6.4})).
The sake for such a condition is due to
estimate (\ref{7.25}) where
a factor $\rho^2$ appear
in front of $\|u\|_{L_{p,-1}(G)}$. Since this factor takes origin
from the definition of $J_7(u, {\cal E}, (x_0,x))$
in (\ref{6.14}), we can say that
the restriction to considering small $\ve$
is a direct consequence
of the necessity of introducing  the
boundary operator ${\cal B}((x_0,x);D_x,D_{x_0})$
in order to prove Lemma \ref{lem4.2}.}
\end{remark}
\section{Appendix}
\setcounter{equation}{0}
We recall here the condition of ellipticity in the sense of
Agranovich-Vishik for the operator
${\cal U}(0,\om,z, D_{\om})=\{L(0,\om,z, D_{\om}),B(0,\om,z, D_{\om})\}$,
$\om\in\Om$, $z\in\csp$,
introduced in Section 3.
Moreover, taking advantage from the discussions on pages 88--90 in
\cite{AV}, we sketch out how easily problems
satisfying this condition can be constructed.\\
First of all, let $G$ be a bounded domain of $\rsp^n$ whose boundary
$\partial G$ is an $(n-1)$--dimensional smooth surface
locally admitting rectification
by means of a $C^{\infty}$ transformation
of co-ordinates $x\to y$. As a result of such transformation
$\partial G$ becomes locally a hyperplane with equation $y_n=0$
and $G$ turns out to lie in the half-space $y_n>0$.
\med\pn
Suppose now we are given the boundary value problem
\beqn\label{8.1}
&&\hskip -1truecm
{\cal L}(x;D_x,q)u(x,q)\,=\,f(x,q)\,,\qq x\in G\,,
\\[2mm]
\label{8.2}
&&\hskip -1truecm
{\cal B}(x'; D_x,q)u(x',q)\,=\,g(x',q)\,,
\qq x'\in\partial G\,.
\eeqn
Here ${\cal L}$ and ${\cal B}$ are
$k\times k$ and $m\times k$
matrix differential operators
with sufficiently smooth complex coefficients
polynomially depending on a parameter $q$ which varies in the sector
$Q=\{z\in \csp:\theta_0\le\arg z\le \theta_1\}$. In particular,
for $\theta_0=\theta_1$, $Q$ can be a ray.\\
The assumptions on the orders of the differential operators
being the same as those in Section 3, with
${\cal L}_0$ and ${\cal B}_0$ we denote here the principal
parts consisting of the terms of higher order
in ${\cal L}$ and ${\cal B}$, respectively.
We impose two algebraic conditions
on ${\cal L}$ and ${\cal B}$.
\vskip 0,3truecm
\indent
{\it i)}\;
If $x\in\ov G$, $\xi\in\rsp^n$ and $q\in Q$,
$|\xi|+|q|\neq 0$, then
\beqn\no
&&\hskip -3truecm
\det {\cal L}_0(x;\xi,q)\neq 0\,.
\eeqn
Since the degree of the polynomial $\det {\cal L}_0(x;\xi,q)$
in $\xi$ is $2r$, for $n\ge 2$
the equation $\l\to\det {\cal L}_0(x;\xi+\l\xi_0,q)=0$,
where $\xi_0\neq 0$ and $\xi$ is orthogonal to $\xi_0$,
has exactly $r$ roots with positive imaginary part.
\vskip 0,3truecm
\indent
{\it ii)}\;
Let $x'$ be any point on $\partial G$.
We consider the problem in the half-line
\beqn
\label{8.3}
&&\hskip -1truecm
{\cal L}_0((x',0);\xi',-iD_y,q)v(y)\,=0\,,\qq y=x_n>0\,,
\\[2mm]
\label{8.4}
&&\hskip -1truecm
{\cal B}_0((x',0);\xi',-iD_y,q)v(y)|_{y=0}\,=\,h\,.
\eeqn
and we require that if $|\xi'|+|q|\neq 0$, $q\in Q$, for any
vector $h\in\csp^k$
this problem has one and only one solution in the class ${\mathfrak M}(\xi')$
of stable solutions of (\ref{8.3}), i.e. solutions tending (exponentially)
to zero together with all their derivatives as $y\to+\infty$.
\med
\pn
For $q=0$ conditions {\it i)} and {\it ii)} reduces, respectively,
to the condition that system (\ref{8.1}) is elliptic and
to the condition of Shapiro-Lopatinskij for problem (\ref{8.1}), (\ref{8.2}).
\begin{definition}
\emph{With the problem (\ref{8.1}), (\ref{8.2}) we associate the operator
\beqn
&&\hskip -1truecm
{\cal U}(x;D_x,q)
=\{{\cal L}(x;D_x,q), {\cal B}(x';D_x,q)\}\,.\no
\eeqn
If ${\cal L}$ and ${\cal B}$ satisfy the algebraic
conditions {\it i)} and {\it ii)},
we say that ${\cal U}(x;D_x,q)$ is elliptic with parameter
in the sense of Agranovich-Vishik.}
\end{definition}
\pn
Examples of problems satisfying {\it i)} and {\it ii)}
can be constructed as follows.\\
Let ${\cal L}(x;D_x,D_{x_{n+1}})$ be an elliptic operator
in the closure of the infinite cylinder
$G_1=G\times( -\infty,+\infty)$,
connected on
$\partial\mbox{}G_1=\partial G\times( -\infty,+\infty)$
by the condition of Shapiro-Lopatinskij
with the boundary operator
${\cal B}(x;D_x,D_{x_{n+1}})$.
Then, the operators ${\cal L}(x;D_x,q)$ and ${\cal B}(x;D_x,q)$, obtained
by replacing $D_{x_{n+1}}$ with $q$, satisfy conditions {\it i)}
and {\it ii)} in each section $x_{n+1}={\rm const}$ of $G_1$
if $q$ belongs to $\{z\in \csp: \arg z=0\}$
or $\{z\in \csp: \arg z=\pi\}$.\\
Moreover, since it is well-known (cf. \cite{ADN2} or \cite{AES})
that the Shapiro-Lopatinskij condition is equivalent
to require that ${\cal B}$ cover ${\cal L}$
on $\partial G$ in the sense of \cite{ADN2}
the former example shows that the uniform ellipticity
of system (\ref{3.16})
implies the ellipticity in the sense of Agranovich--Vishik
for the operator
${\cal U}(0,\om,z, D_{\om})$
defined through (\ref{3.17}), (\ref{3.18}).

\end{document}